\documentclass[11pt]{article}
\usepackage{amsmath,amsthm,amssymb,amsfonts}
\usepackage{multicol, latexsym}
\usepackage{latexsym}
\usepackage{psfrag,import}
\usepackage[final]{graphicx}
\usepackage{fullpage}
\usepackage{framed}
\usepackage{verbatim}
\usepackage{color}
\usepackage{epsfig}
\usepackage[outdir=./]{epstopdf}
\usepackage{hyperref}   
\usepackage{geometry}
\usepackage{mathtools}
\usepackage{enumerate}
\usepackage{enumitem}   
\usepackage{multicol}
\usepackage{a4wide}
\usepackage{booktabs}
\usepackage{enumitem}
\usepackage{lineno}
\usepackage{parcolumns}
\usepackage{thmtools}
\usepackage{xr}
\usepackage{epstopdf}
\usepackage{mathrsfs}
 \usepackage{subcaption}
 \usepackage{soul}
\usepackage{comment}
\usepackage{authblk}
\usepackage{setspace}

\theoremstyle{plain}
\newtheorem{theorem}{Theorem}[section]

\newtheorem{test}[theorem]{Test}

\theoremstyle{definition}
\newtheorem{definition}[theorem]{Definition}

\theoremstyle{remark}

\newcommand\SP{\mathcal{S}}

\newcommand\R{\mathbb{R}}

\newcommand\bx{\boldsymbol{x}}
\newcommand\bw{\boldsymbol{w}}
\newcommand\bk{\boldsymbol{k}}
\newcommand\bp{\boldsymbol{p}}
\newcommand\by{\boldsymbol{y}}

\newcommand\bz{\boldsymbol{z}}

\newcommand\bu{\boldsymbol{u}}

\graphicspath{images/}
        {\begin{list}
                {\noindent\makebox[0mm][r]{$\bullet$}}
                {\leftmargin=5.5ex \usecounter{enumi}
      \topsep=1.5mm \itemsep=-.75ex}
        }
        {\end{list}}
        
\title{Endotactic and strongly endotactic networks with infinitely many positive steady states}

\author[1]{Samay Kothari \thanks{samay.kothari@research.iiit.ac.in}}
\author[2]{Abhishek Deshpande \thanks{abhishek.deshpande@iiit.ac.in}}
\affil[1,2]{Center for Computational Natural Sciences and Bioinformatics, International Institute of Information Technology Hyderabad}

\begin{document}

\maketitle

\begin{abstract}

%\textcolor{red}{The number of steady states of a reaction network is an important metric while studying its behaviour.} 
We show that there exists endotactic and strongly endotactic dynamical systems that are not weakly reversible and possess infinitely many steady states. We provide a few examples in two dimensions and an example in three dimensions that satisfy this property. In addition, we prove for some of these systems that there exist no weakly reversible mass-action systems that are dynamically equivalent to mass-action systems generated by these networks. 

\end{abstract}

\section{Introduction}\label{sec:intro}

Polynomial dynamical systems have profound applications in the study of ecology, epidemiology, and population dynamics. These dynamical systems are often of the form 

\begin{eqnarray}
\frac{d\bx}{dt} = \displaystyle\sum_{i=1}^m{\bx}^{\by_i}\bw_i
\end{eqnarray}
where $\bx\in\mathbb{R}^n_{>0}, \by_i,\bw_i \in\mathbb{R}^n$ and ${\bx}^{\by_i} = x_1^{y_1}x_2^{y_2},..., x_m^{y_m}$. In particular, properties like extinction of species, persistence (property that no species goes extinct) and permanence (convergence to a compact set) are of paramount importance~\cite{craciun2013persistence,gopalkrishnan2014geometric,craciun2015toric}. Many of these dynamical properties arise as a manifestation of the underlying structure of the reaction network. For example, it is conjectured that weakly reversible reaction networks (i.e., networks where every reaction is part of a directed cycle) are permanent. This is known as the \emph{Permanence} conjecture~\cite{craciun2013persistence}. Another important component in the analysis of dynamical systems is the study of steady states. We outline some important properties of steady states below.

Questions like the number of positive steady states, their stability and bifurcation are crucial to their analysis. In this paper, our focus is on the number of positive steady states. The analysis of steady states has a long history starting from the work on Deng et. al.~\cite{deng2011steady}, where they studied the existence and finiteness of steady states for weakly reversible networks. It has now been established by Boros~\cite{boros2019existence} that for weakly reversible networks, there exists a positive steady state within each stoichiometric compatibility class. Further Boros, Craciun and Yu~\cite{boros2020weakly} have shown that there exists weakly reversible mass-action systems that possess infinitely many positive steady states; in particular they have constructed examples of weakly reversible mass-action systems where there is a curve of positive steady states. In this paper, we show that there exists examples of endotactic and strongly endotactic networks that possess a curve of positive steady states. 
%\textcolor{red}{The endotactic reaction networks have lesser constraints to satisfy as compared to the weakly-reversible reaction networks, and strongly endotactic networks show other properties such as meta-stability, emergent behaviour and non-linear dynamics.} 
The networks that appear in our construction are not weakly reversible; and in certain cases there exists no weakly reversible mass-action systems that are dynamically equivalent to the mass-action systems generated by these networks.  In addition, we construct three dimensional endotactic networks that are not weakly reversible and possess infinitely many positive steady states. 

In Section~\ref{sec:reaction_networks} we define reaction networks as graphs embedded in Euclidean space. In addition, we review some concepts from reaction network theory. In Section~\ref{sec:steady_states} we review certain properties about steady-states for weakly reversible and strongly endotactic networks. In Section~\ref{sec:dynamical_equivalence} we define \emph{dynamically equivalent networks}. In Section~\ref{sec:operations} we illustrate a set of simple operations on reaction networks that can be used to transform them to more complicated systems. In Section~\ref{sec:strongly_endotactic} we give examples of strongly endotactic and endotactic networks that are not weakly reversible but possess infinitely many steady states. In particular, we describe operations on a \emph{base unit network} that eventually generate these endotactic and strongly endotactic networks. In Section~\ref{sec:discussion} we review our results and chalk out directions for future research.

\section{Reaction networks}\label{sec:reaction_networks}

A reaction network is a directed graph where reactions are represented by edges; reactant and product complexes are represented by vertices. Such graphs have been referred to as Euclidean embedded graphs(or E-graphs)~\cite{craciun2015toric,craciun2019polynomial,craciun2020endotactic,craciun2019quasi} in literature. We will denote this graph by $G=(V,E)$, where $V$ is the set of vertices and $E$ is the set of edges. We will assume that $V\subset\R^n$ and $E\subset V\times V$. We will denote the set of species by $\SP$. A reaction $\by\rightarrow\by'$ will be denoted by the edge $(\by,\by')\in E$. We will denote by $V_{\rm{source}}$ the set of source vertices. We define the following terms with respect to the reaction network $G=(V,E)$:

\begin{enumerate}

\item The \emph{stoichiometric subspace} of a reaction network is the set $S:=\{\rm{span}(\by'-\by)\, |\, \by\rightarrow\by'\in E\}$. We will denote by $s$ the dimension of the stoichiometric subspace $S$.

\item The \emph{stoichiometric compatibility class} of a point $\bz_0\in\R^n_{>0}$ is the set $S_{\bz_0}:=\{(\bz_0 + S)\cap\R^n_{>0}\}$.

\item A set $V_0\subseteq V$ is called a \emph{linkage class} if $V_0$ is a maximal connected component. We will denote the number of linkage classes in a network by $\ell$.

\item The deficiency of a network(denoted by $\delta)$ is given by the formula: $\delta = |V| - \ell - s$. Note that $\delta$ is an integer such that $\delta\geq 0$~\cite{gunawardena2003chemical}. 

\item A reaction network is \emph{weakly reversible} if each reaction in the network is part of a directed cycle.

\item A reaction network is \emph{endotactic}~\cite{craciun2013persistence,anderson2020classes} if for every reaction $\by\rightarrow\by'\in E$ and $\bu\in\R^n$ that satisfies $\bu\cdot(\by'-\by)<0$, there exists $\bar{\by}\rightarrow\bar{\by}'\in E$ such that $\bu\cdot\bar{\by} < \bu\cdot\by$ and $\bu\cdot(\bar{\by}'- \bar{\by})>0$.

\item A reaction network is \emph{strongly endotactic}~\cite{gopalkrishnan2014geometric,anderson2020classes} if for every reaction $\by\rightarrow\by'\in E$ and $\bu\in\R^n$ that satisfies $\bu\cdot(\by'-\by)<0$, there exists $\bar{\by}\rightarrow\bar{\by}'\in E$ such that $\bu\cdot\bar{\by} < \bu\cdot\by$, $\bu\cdot(\bar{\by}'- \bar{\by})>0$ and $\bu\cdot\bar{\by}\leq \bu\cdot\tilde{\by}$ for every $\tilde{\by}\in V_{\rm{source}}$.

\item A set of species $\SP'\subseteq \SP$ is called a \emph{siphon} if for every reaction $\by\rightarrow\by'$ the following holds true: 
if the complex $\by$ does not contain any species from $\SP'$, then the complex $\by'$ also does not contain any species from $\SP'$. 

\item A set of species $\SP'\subseteq \SP$ is said to be \emph{critical} if there exists a point $\bp\in\R^n_{\geq 0}$ with $\bp_i=0$ if $i\in\SP'$ such that $(\bp + S)\cap\R^n_{>0}\neq\emptyset$.

\end{enumerate} 

Figure~\ref{fig:E_graph} lists several examples of reaction networks that illustrates these properties.

\begin{figure}[h!]
\centering
\includegraphics[scale=0.45]{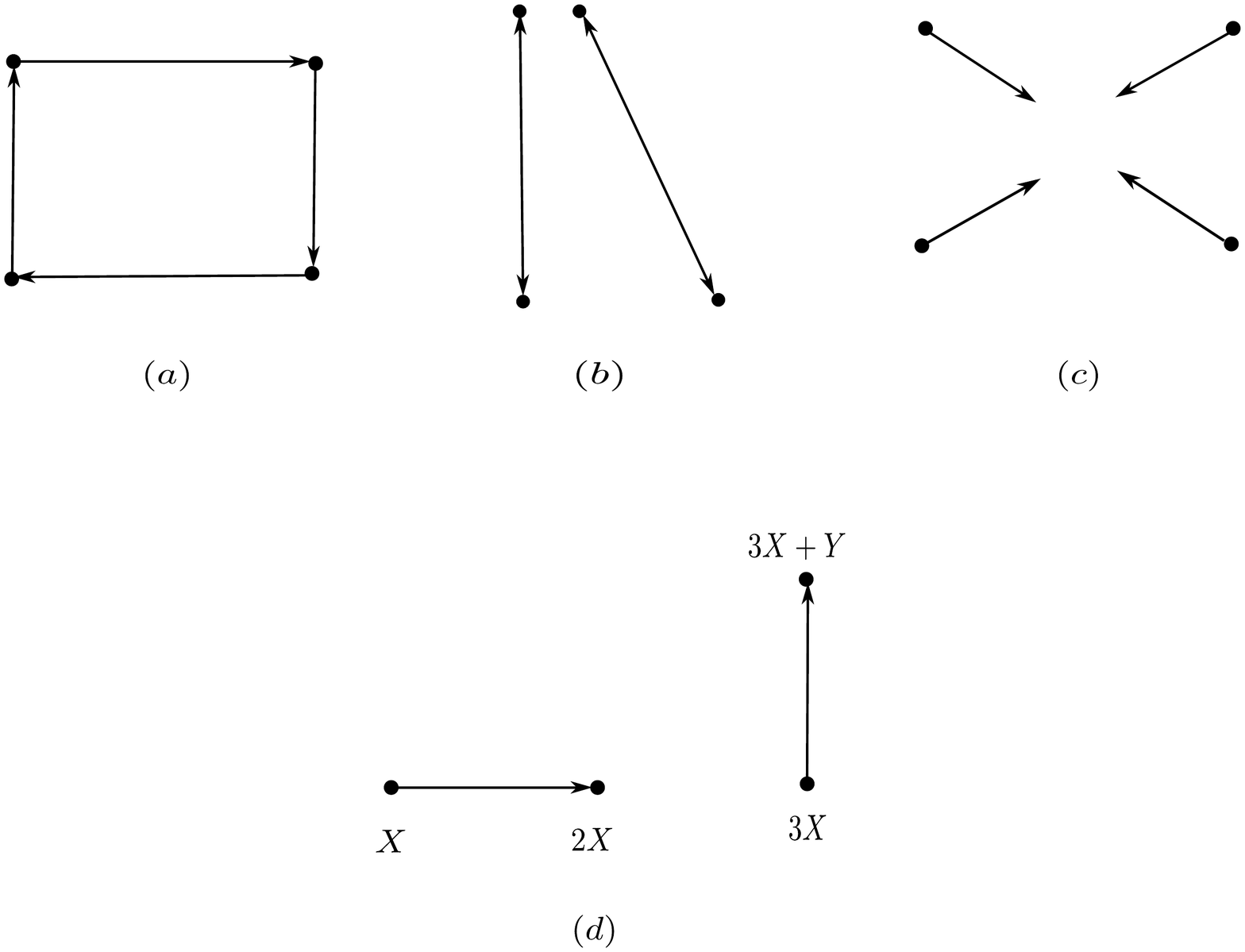}
\caption{\small (a) A weakly reversible reaction network. It consists of a single linkage class and has a deficiency $\delta=1$. (b) This is an endotactic network that is not strongly endotactic. It possesses two linkage classes and has a deficiency $\delta=0$. (c) This is an example of a strongly endotactic network. It possesses four linkage classes and has a deficiency $\delta=2$. (d) This network possesses two linkage classes and has a deficiency $\delta=0$. This network has the following sets as siphons: $\{X\},\{X,Y\}$; the set $\{Y\}$ is not a siphon. The sets $\{X\},\{Y\},\{X,Y\}$ are critical.}
\label{fig:E_graph}
\end{figure}

The notions of endotactic and strongly endotactic network can be understood using the \emph{parallel sweep test}~\cite{craciun2013persistence,gopalkrishnan2014geometric}. We state the test below for completeness.

\begin{test}[Parallel sweep test for strongly endotactic networks]\label{test_strongly_endotactic}
Consider a reaction network $G=(V,E)$. Let $S_G$ denote the stoichiometric subspace and $\rm{conv}(V_{\rm{source}})$ denote the convex hull of the source vertices of $G$. Choose a vector $\bu\in\mathbb{R}^n$ satisfying $\bu\notin {S_G}^{\perp}$. Now sweep with a hyperplane perpendicular to $\bu$ towards $\rm{conv}(V_{\rm{source}})$. Let $H_0$ denote the hyperplane when it first touches the boundary of  $\rm{conv}(V_{\rm{source}})$. Then the network is said to pass the parallel sweep for $\bu$ if the following conditions hold:
\begin{enumerate}
\item For all reactions $\by\rightarrow \by'\in E$ with source $\by\in H_0$, we have $\by\cdot (\by'-\by) \geq 0$ and
\item  There exists at least one reaction $\hat{\by}\rightarrow \hat{\by}'\in E$ with $\hat{\by}\in H_0$ such that $\bu\cdot (\hat{\by}'- \hat{\by}) >0$
\end{enumerate}
We say that the network is strongly endotactic if it passes the parallel sweep testfor all $\bu$. Otherwise, the network fails the parallel sweep test and is not strongly endotactic.
\end{test}

\begin{figure}[h!]
\centering
\includegraphics[scale=0.35]{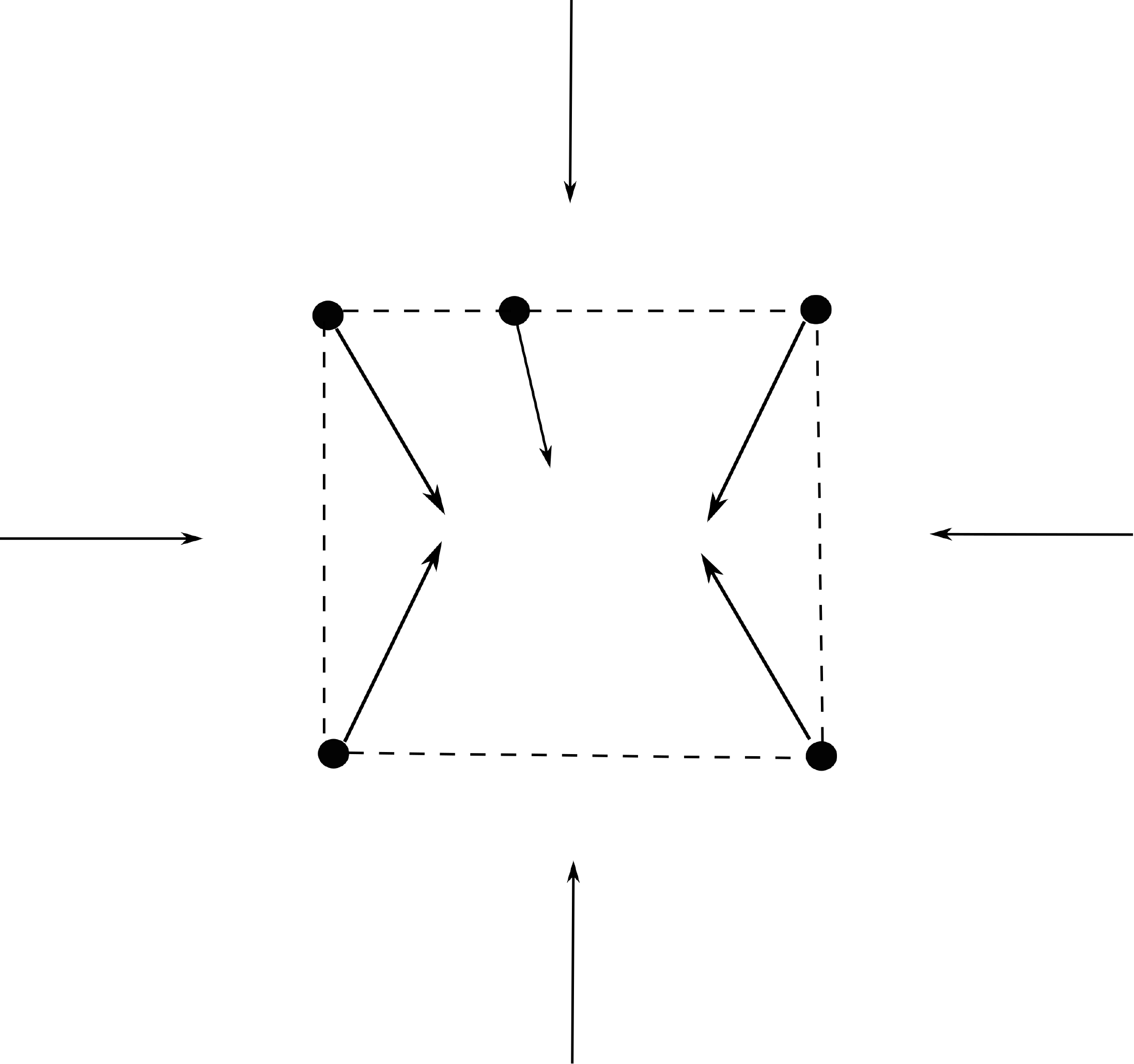}
\caption{\small The boundary of the convex hull of this network is shown by dotted lines. This is an example of a strongly endotactic network since it passes the parallel sweep test~\ref{test_strongly_endotactic} . A few candidate vectors $\bu$ are shown in the figure.}
\label{fig:strongly_endotactic}
\end{figure}

\begin{test}[Parallel sweep test for endotactic networks]\label{test_endotactic}
Consider a reaction network $G=(V,E)$. Let $S_G$ denote the stoichiometric subspace and $\rm{conv}(V_{\rm{source}})$ denote the convex hull of the source vertices of $G$. Choose a vector $\bu\in\mathbb{R}^n$. Now sweep with a hyperplane perpendicular to $\bu$ towards $\rm{conv}(V_{\rm{source}})$. Let $H_0$ denote the hyperplane when it first touches the boundary of  $\rm{conv}(V_{\rm{source}})$. We first verify if the following  holds: for all reactions $\by\rightarrow \by'\in E$ with source $\by\in H_0$, we have $\by\cdot (\by'-\by) \geq 0$. If no, then the network fails the parallel sweep test and is not endotactic. If yes, then we have two cases:

\begin{enumerate}
\item If there exists at least one reaction $\hat{\by}\rightarrow \hat{\by}'\in E$ with $\hat{\by}\in H_0$ such that $\bu\cdot (\hat{\by}'- \hat{\by}) >0$, then the network passes the parallel sweep test for $\bu$.

\item For all reactions $\by\rightarrow \by'\in E$ with source $\by\in H_0$, if we have $\bu\cdot (\by'-\by) =0$, then we continue sweeping withe the hyperplane and repeat the steps above.  
\end{enumerate}
We say that the network is endotactic if it passes the parallel sweep test for all $\bu\in\mathbb{R}^n$.
\end{test}

\begin{figure}[h!]
\centering
\includegraphics[scale=0.45]{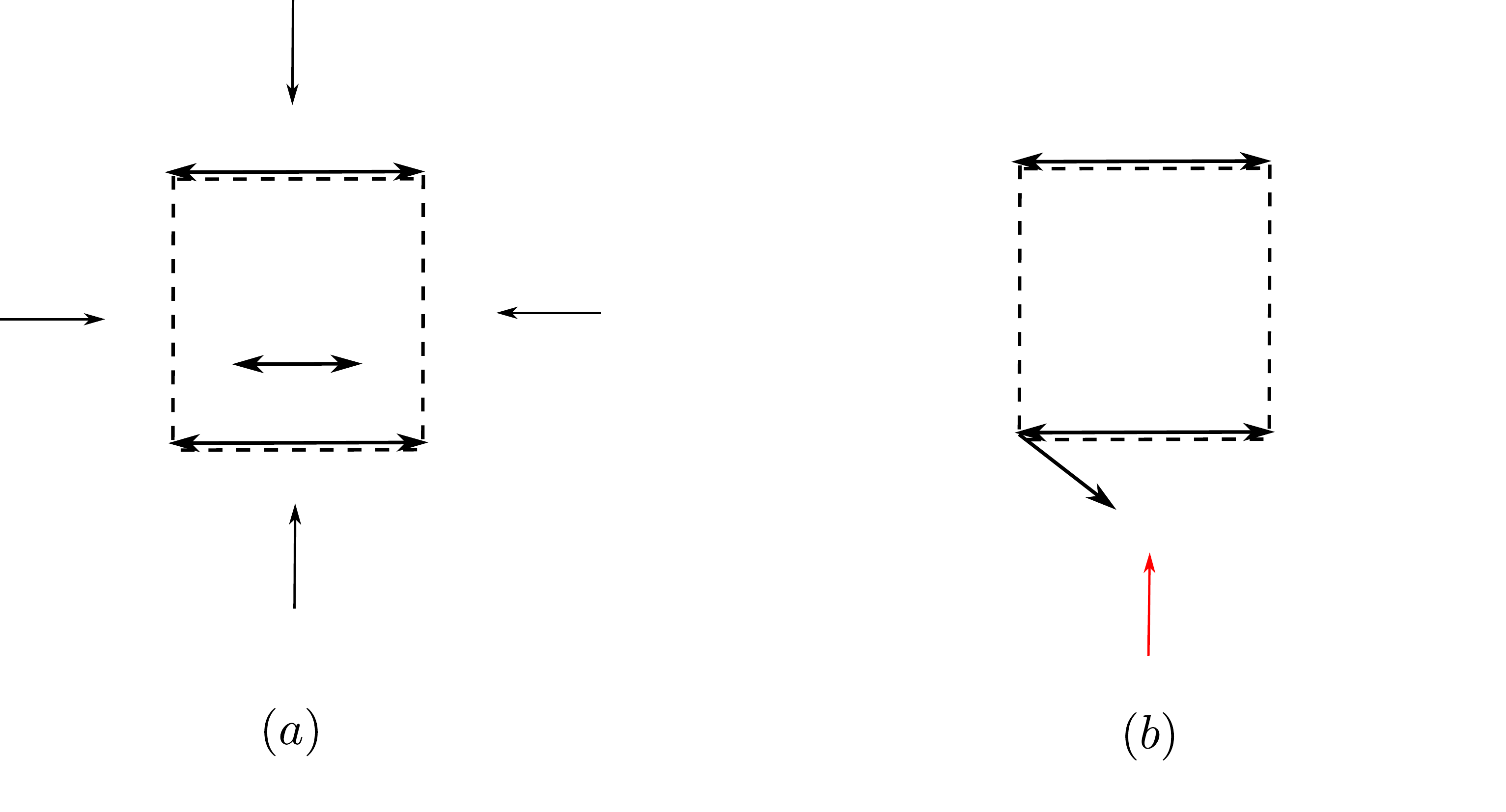}
\caption{\small (a) The boundary of the convex hull of this network is shown by dotted lines. This is an example of endotactic network that is not strongly endotactic. It passes test~\ref{test_endotactic}. A few candidate vectors $\bu$ are shown in the figure. (b) The boundary of the convex hull of this network is shown by dotted lines. This network is not endotactic since it fails the parallel sweep test~\ref{test_endotactic} when we sweep in the direction marked in red (i.e., when $\bu$ is the direction marked in red).}
\label{fig:endotactic}
\end{figure}

Note that weakly reversible reaction networks are endotactic. Strongly endotactic reaction networks are also endotactic. Further weakly reversible reaction networks consisting of a single linkage class are strongly endotactic. However, there exists endotactic reaction networks  that are not weakly reversible; and endotactic networks that are not strongly endotactic.
%\textcolor{red}{Note that while weakly reversible reaction networks are considered endotactic, strongly endotactic reaction networks are also classified as endotactic. Also a weakly reversible reaction network comprising a single linkage class is classified as strongly endotactic. However, it is important to acknowledge that there exist endotactic reaction networks that do not fall under the category of weakly reversible or strongly endotactic reaction networks, as well strongly endotactic networks that are not weakly reversible.}

Brunner et.al.~\cite{anderson2020classes} define \emph{extremally weakly reversible} networks as those networks whose projections relative to the boundary of the convex hull of the source vertices are weakly reversible. They show that strongly endotactic networks in two dimensions are dynamically equivalent to some extremally weakly reversible network. 

Reaction networks can be assigned a set of rate constants (one rate constant per reaction) which regulates the rate at which the species are produced or consumed. Given as assignment of rates, each reaction network can be linked to a set of dynamical system. In particular, if the dynamics is governed by mass-action kinetics~\cite{feinberg1979lectures,voit2015150,guldberg1864studies,yu2018mathematical,gunawardena2003chemical,adleman2014mathematics},  then the rate equations are given by the following:

\begin{eqnarray}\label{eq:mass_action}
\frac{d\bx}{dt} = \displaystyle\sum_{\by\rightarrow \by'\in E}{\bk}_{\by\rightarrow \by'}{\bx}^{\by}(\by'- \by)
\end{eqnarray}
A point $\bar{\bx}\in\R^n_{\geq 0}$ is called a steady state of~\ref{eq:mass_action} if $\displaystyle\sum_{\by\rightarrow \by'\in E}{\bk}_{\by\rightarrow \by'}{\bar{\bx}}^{\by}(\by'- \by)=\textbf{0}$. We will denote this mass-action dynamical system by $G_{\bk}$. A mass-action system $G_{\bk}$ is said to be \emph{complex balanced} if  there exists $\tilde{\bx}\in\mathbb{R}^n_{\geq 0}$ such that the following is true for every vertex $\by_0$:
\begin{eqnarray}
\displaystyle\sum_{\by_0\rightarrow \by'\in E}{\bk}_{\by_0\rightarrow \by'}{\tilde{\bx}}^{\by_0} = \displaystyle\sum_{\by'\rightarrow \by_0\in E}{\bk}_{\by'\rightarrow \by_0}{\tilde{\bx}}^{\by'}
\end{eqnarray}
A lot of work hovers around inferring the dynamics of a reaction network from its underlying graphical structure. Such dynamical properties include \emph{persistence}, \emph{permanence} and \emph{global stability}. The idea of persistence means that for any solution $\bx(t)$ with $\bx(0)\in\R^n_{>0}$, we have $\displaystyle\lim\inf_{t\to\infty}\bx_i(t)>0$ for all $i\in\SP$. A dynamical system is \emph{structurally persistent} if the underlying reaction network does not possess critical siphons~\cite{angeli2007petri_2,ali2020computational}. Permanence is a stronger condition than persistence; in particular it implies that for every stoichiometric compatibility class $D$, there exists a compact set $J\subset D$ such that the solution satisfies $\bx(t)\in J$ for large enough $t$. A point $\bx^*\in\R^n_{\geq 0}$ is a global attractor if for any solution $\bx(t)$ with $\bx(0)\in\R^n_{>0}$, we have $\displaystyle\lim_{t\to\infty}\bx(t) = \bx^*$.

It is conjectured that weakly reversible reaction networks are permanent~\cite{craciun2013persistence}. Further it is also conjectured that complex balanced systems have a globally attracting steady state within each stoichiometric compatibility class. This is known as the \emph{Global Attractor Conjecture}. In recent years, there has been considerable progress in proving these conjectures. In particular, Craciun, Nazarov and Pantea~\cite{craciun2013persistence} have proved that endotactic networks in two dimensions are permanent. This has been extended by Pantea~\cite{pantea2012persistence} to networks with two dimensional stoichiometric subspace. Anderson~\cite{anderson2011proof} has shown that weakly reversible networks with a single linkage class are persistent. Gopalkrishnan, Miller and Shiu have shown that strongly endotactic networks are permanent. Angeli and Sontag~\cite{angeli2007petri} have shown that weakly reversible reversible networks that do not possess critical siphons are persistent. The entire proof of the Global Attractor Conjecture has been proposed in 2015 by Craciun~\cite{craciun2015toric}.

\section{Steady states}\label{sec:steady_states}

%\textcolor{red}{A steady state is a state in which the concentrations of all chemical species remain constant over time, i.e., rate of change of concentrations for all the species is zero.}
To understand the dynamics exhibited by reaction networks, it is important to analyze the behaviour of their steady states. In particular, phenomena like the existence of a unique steady state, or the the existence of multiple steady states or oscillations and bifurcations are of special interest. Proving the existence of positive steady states has received a lot of attention from the reaction network community. In particular, it is known that complex balanced dynamical systems have a unique steady state in each of their stoichiometric compatibility classes~\cite{horn1972general}. Further, these steady states are locally asymptotically stable. Boros~\cite{boros2019existence} has shown that for weakly reversible networks, there exists a steady state within each stoichiometric compatibility class. Gopalkrishnan et.al.~\cite{gopalkrishnan2014geometric} have shown that strongly endotactic networks possess a positive steady state within each stoichiometric compatibility class.

In recent work, Craciun, Boros and Yu~\cite{boros2020weakly} have shown that there exists weakly reversible networks that possess infinitely many steady states. In particular, they have shown examples of weakly reversible networks that possess a curve of steady states. Our work extends this to show that there exists endotactic and strongly endotactic networks which are not weakly reversible, but possess a curve of steady states. Examples in Section~\ref{sec:steady_states} illustrate this point.

\section{Dynamical Equivalence}\label{sec:dynamical_equivalence}

The phenomenon of \emph{dynamical equivalence} has been studied extensively in recent years. It has been also called \emph{macroequivalence} by Horn and Jackson~\cite{horn1972general} and \emph{confoundability} by Craciun and Pantea~\cite{craciun2008identifiability}. We formally define it below.

\begin{definition}
Given dynamical systems $G_{\bk}$ and $G'_{\bk'}$, we say that they are dynamically equivalent if the following holds for every $\bx\in\mathbb{R}^n_{>0}$
\begin{eqnarray}
\displaystyle\sum_{\by\rightarrow \by'\in E}{\bk}_{\by\rightarrow \by'}{\bx}^{\by}(\by'- \by) = \displaystyle\sum_{\bar{\by}\rightarrow \bar{\by}'\in E'}{\bk'}_{\bar{\by}\rightarrow \bar{\by}'}{\bx}^{\bar{\by}}(\bar{\by}'-\bar{\by})
\end{eqnarray}
\end{definition}

\begin{figure}[h!]
\centering
\includegraphics[scale=0.45]{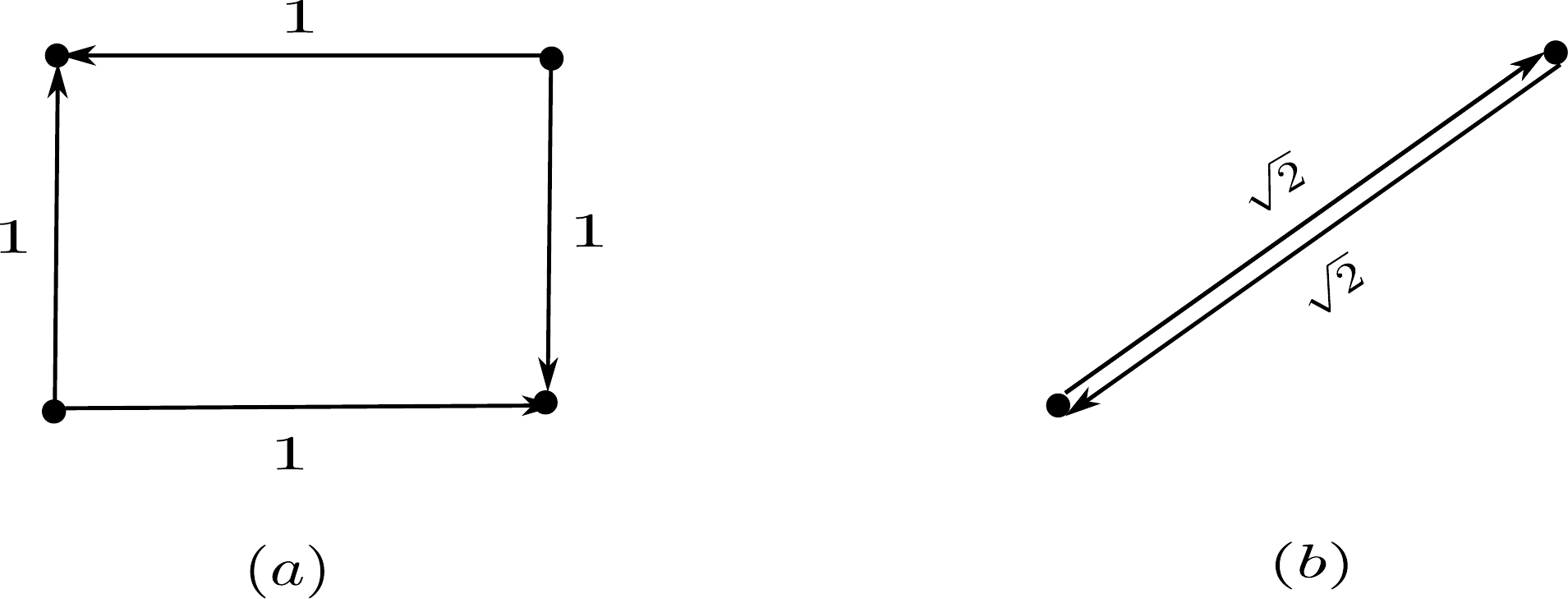}
\caption{\small The mass-action system generated by the network in (a) is dynamically equivalent to the mass-action system generated by the network in (b).}
\label{fig:E_graph}
\end{figure}

If two dynamical systems  $G_{\bk}$ and $G'_{\bk'}$ are dynamically equivalent, we will denote this by $G_{\bk}\sim G'_{\bk'}$.

\section{Operations on Mass Action Systems}\label{sec:operations}

In this section, we describe different operations that we would be using on the mass-action systems in Section~\ref{sec:strongly_endotactic} such as: translation, scalar multiplication and addition.

\subsection{Translation}

For a mass-action system, when we multiply the right hand side of the dynamical system with by a monomial $x^ay^b$, it translates the system by a vector $(a,b)^T$ on the graph. For example consider the dynamical system given by:

\begin{figure}
    \centering
    \includegraphics[width = 0.2\textwidth]{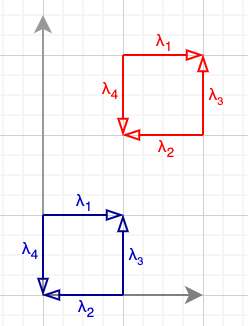}
    \caption{Translations of a mass action system}
    \label{fig:translations}
\end{figure}
\begin{equation}\label{translation}
     \begin{aligned}
        \dot{x}&=(\lambda_1y-\lambda_2x) \\
        \dot{y}&=(\lambda_3x-\lambda_4y)
    \end{aligned}
\end{equation}

When we multiple the right-hand side of Equation~\ref{translation} is multiplied by monomial $xy^2$, it translates the system by vector $(1,2)^T$, and the dynamical system is given by:

\begin{equation*}
     \begin{aligned}
        \dot{x}&=xy^2(\lambda_1y-\lambda_2x) \\
        \dot{y}&=xy^2(\lambda_3x-\lambda_4y)
    \end{aligned}
\end{equation*}
Figure \ref{fig:translations} represents the translation operation by a vector $(1,2)^T$.

\subsection{Scalar Multiplication}

The scalar multiplication of the rate constant of the mass action system by a factor of $\rho$ give rise to two different behaviours depending upon the sign of the $\rho$.

When $\rho >0$, we modify the rate constants of the system by multiplying each of them with $\rho$. For example, let us assume we have an initial mass action system as:

\begin{equation}\label{scalar}
    \begin{aligned}
        \dot{x}&=(\lambda_1xy^2-\lambda_2x^2y) \\
        \dot{y}&=(\lambda_3x^2y-\lambda_4xy^2)
    \end{aligned}
\end{equation}

\begin{figure}[h!]
    \centering
    \begin{subfigure}[t]{0.33\textwidth}
        \centering
        \includegraphics[width = 0.7\textwidth]{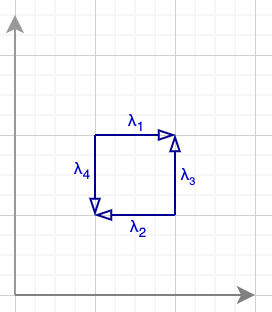}
        \caption{Initial System}
        \label{fig:scalarInitialSystem}
    \end{subfigure}%
    \begin{subfigure}[t]{0.33\textwidth}
        \centering
        \includegraphics[width = 0.7\textwidth]{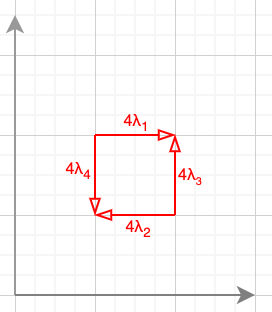}
        \caption{$\rho = 4$}
        \label{fig:scalarMultiplication1}
    \end{subfigure}%
    \begin{subfigure}[t]{0.33\textwidth}
        \centering
        \includegraphics[width = 0.7\textwidth]{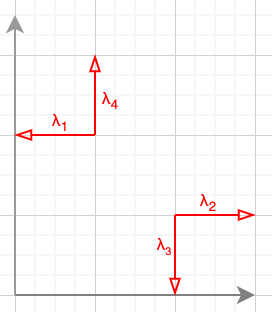}
        \caption{$\rho = -1$}
        \label{fig:scalarMultiplication2}
    \end{subfigure}%
    \caption{Scalar Multiplication of a mass-action system}
\end{figure}

We then multiply this dynamical system by $\rho = 4$ on the right-hand side of Equation~\ref{scalar}, we get the following dynamical system
\begin{equation*}
    \begin{aligned}
        \dot{x}&=4(\lambda_1xy^2-\lambda_2x^2y) \\
        \dot{y}&=4(\lambda_3x^2y-\lambda_4xy^2)
    \end{aligned}
\end{equation*}

Figure~\ref{fig:scalarMultiplication1} shows the dynamical system obtained by multiplying the system given by Equation~\ref{scalar} by $\rho =4$.

When $\rho <0$, we flip each reaction of the system along it's source vertex and multiply to each rate coefficients by the absolute value of the $\rho$. For example the initial system given by Equation~\ref{scalar} when multiplied by $\rho = -1$ results in a dynamical system given by:
\begin{equation*}
    \begin{aligned}
        \dot{x}&=-1(\lambda_1xy^2-\lambda_2x^2y) \\
        \dot{y}&=-1(\lambda_3x^2y-\lambda_4xy^2)
    \end{aligned}
\end{equation*}

Figure~\ref{fig:scalarMultiplication2} represents the mass-action system obtained by multiplying the system given by Equation~\ref{scalar} by $\rho =-1$.

\subsection{Addition and Simplification}

Addition of two mass-action systems is equivalent to the following: take the union of all the reaction vectors of the two mass-action systems. The rate constants are added for the reaction vectors that are common to both the mass-action systems. For example consider the two mass-action systems:

\begin{equation}\label{sys_1}
    \begin{aligned}
        \dot{x_1}&=(\lambda_1y^3-\lambda_2xy^2) \\
        \dot{y_1}&=(\lambda_3xy^2-\lambda_4y^3) \\
          \end{aligned}
\end{equation}

 \begin{equation}\label{sys_2}
   \begin{aligned}
        \dot{x_2}&=(-\lambda_5xy^2+\lambda_6x^2y) \\
        \dot{y_2}&=(-\lambda_7x^2y+\lambda_8xy^2) \\
    \end{aligned}
\end{equation}

We can also simplify the resultant mass-action system by combining similar monomials. This process will be called \emph{Simplification}. For example, the addition of the mass-action systems given by Equations~\ref{sys_1} and~\ref{sys_2} is given by the following:

\begin{equation*}
    \begin{aligned}
        \dot{x}&=(\lambda_1y^3-(\lambda_2+\lambda_5)xy^2 + \lambda_6 x^2y) \\
        \dot{y}&=((\lambda_3+\lambda_8)xy^2-\lambda_4y^3 - \lambda_7 x^2y) \\
    \end{aligned}
\end{equation*}

\begin{figure}[h!]
    \centering
    \begin{subfigure}[t]{0.33\textwidth}
        \centering
        \includegraphics[width = 0.7\textwidth]{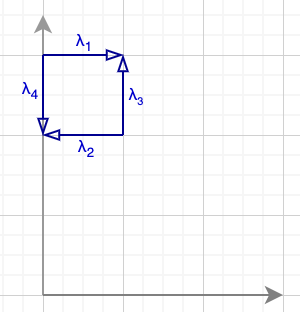}
        \caption{}
        \label{fig:addition1}
    \end{subfigure}%
    \begin{subfigure}[t]{0.33\textwidth}
        \centering
        \includegraphics[width = 0.7\textwidth]{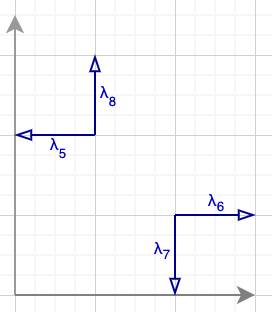}
        \caption{}
        \label{fig:addition2}
    \end{subfigure}%
    \begin{subfigure}[t]{0.33\textwidth}
        \centering
        \includegraphics[width = 0.7\textwidth]{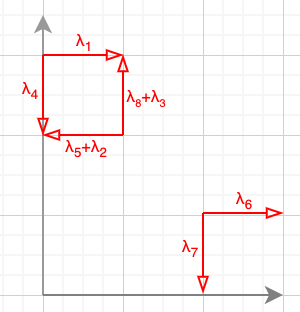}
        \caption{}
        \label{fig:addition3}
    \end{subfigure}%
    \caption{Addition of two mass action systems. By adding system (a) and (b), we get the mass-action system (c)}
\end{figure}

\section{Strongly endotactic and endotactic reaction networks with infinitely many positive steady states}\label{sec:strongly_endotactic}

The goal of this section is to give a few examples of endotactic/strongly endotactic networks with infinitely many positive steady states. The idea of constructing these networks is similar to one described in~\cite{boros2020weakly}. For this process we start with a small unit, a simple endotactic/strongly endotactic network called the \textit{base unit}. We then take multiple base units which are translated by multiplying with different scalar factors, which are combined and \emph{simplified} to form a large reaction network that is endotactic/strongly endotactic but not weakly reversible.

\textit{Simplifying} a reaction network refers to the process of reducing the differential equation to the minimum possible terms in such a way that the resultant reaction network is dynamically equivalent to the mass-action system generated by the original reaction network. The differential equation corresponding to the reaction network, comprises of two components

\begin{itemize}

\item \textbf{Scalar Polynomial}: It refers to common factor that represents the sum of the transformations of the base unit.

% While constructing these polynomials we keep the negative terms in the interior of the convex hull of the positive terms, to construct an object called as a newton polytope. Since we have a single negative term in the polynomial, the value of the polynomial near the boundary of $\mathbb{R}^n_{>0}$ will be positive. If we keep the coefficient of the negative term sufficiently large enough, the value of the polynomial at $(1,1,\ldots,1)^T$ is negative. Consequently, this leads polynomial to have infinitely many roots in $\mathbb{R}^n_{>0}$, thus infinitely many positive steady states.

As remarked in~\cite{boros2020weakly}, the key object here is called the \emph{Newton polytope}~\cite{sturmfels1996grobner}, where the exponent corresponding to the negative monomial lies in the interior of the convex hull of exponents corresponding to the positive monomials. From~\cite{pantea2012global}, this implies that the sign of the polynomial near the boundary is positive. The polynomial evaluated at the point $(1,1,....,1)^T$ can be made negative, by choosing a sufficiently large negative coefficient corresponding to the negative monomial. This change in signs implies that the dynamical system has infinitely many positive steady states.

\item \textbf{Vector Polynomial}: This part of the equation represents the dynamics corresponding to the base unit (which is a simple endotactic/strongly endotactic network) that we are using for construction of a complex reaction network. This base unit is then transformed by the scalar polynomial.
\end{itemize}

\subsection{Example 1}\label{Exp1}

\begin{figure}[h!]
    \centering
    \includegraphics[scale=0.4]{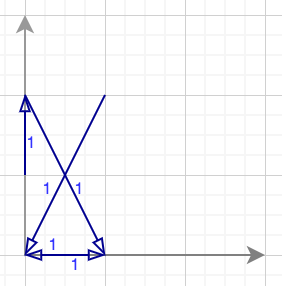}
    \caption{Base Unit 1}
    \label{Base Unit}
\end{figure}

We start with a base unit as given in Figure~\ref{Base Unit}, where the rate constants corresponding to all reactions are set to 1. The network is not weakly reversible, has a deficiency of two and is endotactic, but not strongly endotactic. Further, the network possesses critical siphons implying that the dynamics generated by it is not structurally persistent. The properties mentioned above can be verified with the CoNtRol software~\cite{donnell2014control}.

\begin{figure}[h!]
    \centering
    \includegraphics[width = 0.23\textwidth, height = 150pt]{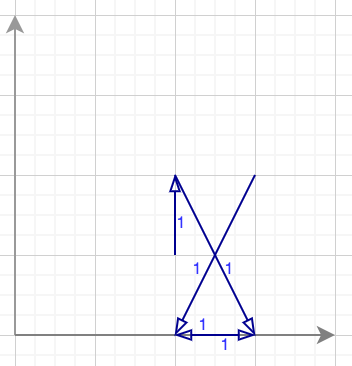}\hfill
    \includegraphics[width = 0.23\textwidth, height = 150pt]{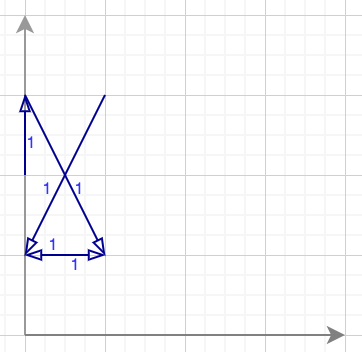}\hfill
    \includegraphics[width = 0.23\textwidth, height = 150pt]{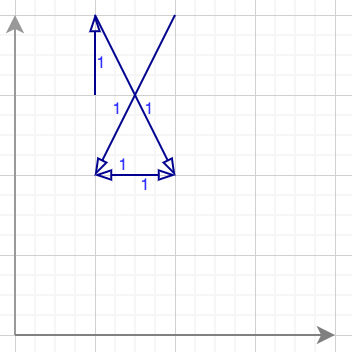}\hfill
    \includegraphics[width = 0.23\textwidth, height = 150pt]{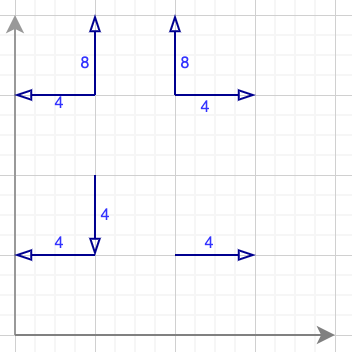}\hfill
    \caption{The different transformations of the base unit given in Figure~\ref{Base Unit}.}
    \label{Transformations}
\end{figure}

\begin{figure}[h!]
    \centering
    \includegraphics[width=0.3\textwidth]{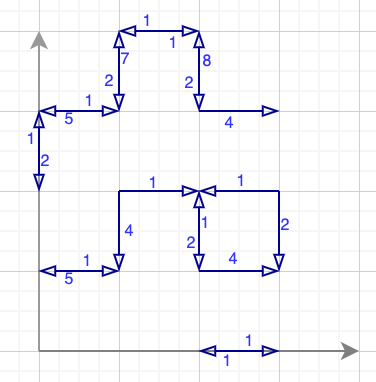}
    \caption{Full unit for Example 1 is obtained by combining all the transformations given in Figure~\ref{Transformations} and simplifying it. The full unit is an endotactic network that is neither strongly endotactic nor weakly reversible.}
    \label{fullSystem1}
\end{figure}

Note that there exists no weakly reversible mass-action system that is dynamically equivalent to the mass-action system generated by the base unit. This base unit is then modified using certain transformations given in Figure~\ref{Transformations}. Combining these transformations, we obtain an endotactic mass action system as depicted in Figure~\ref{fullSystem1}. The dynamics of this system is given by:

\begin{equation*}
    \begin{aligned}
        \dot{x}&=(x^2+xy^2+y-4xy)(1-x+y^2-xy^2) \\
        \dot{y}&=(x^2+xy^2+y-4xy)(y-2y^2-2xy^2)
    \end{aligned}
\end{equation*}
\subsubsection{Phase Plane Analysis}

The positive steady states of the dynamical system are given by:
\begin{equation*}
    \{(x,y) \in \mathbb{R}^2_{>0}:x^2+xy+y-4xy = 0 \}  \cup \{(1,0.25)\}
\end{equation*}

\begin{figure}[h!]
    \centering
    \subfloat[Base Unit]{\includegraphics[width = 0.5\textwidth]{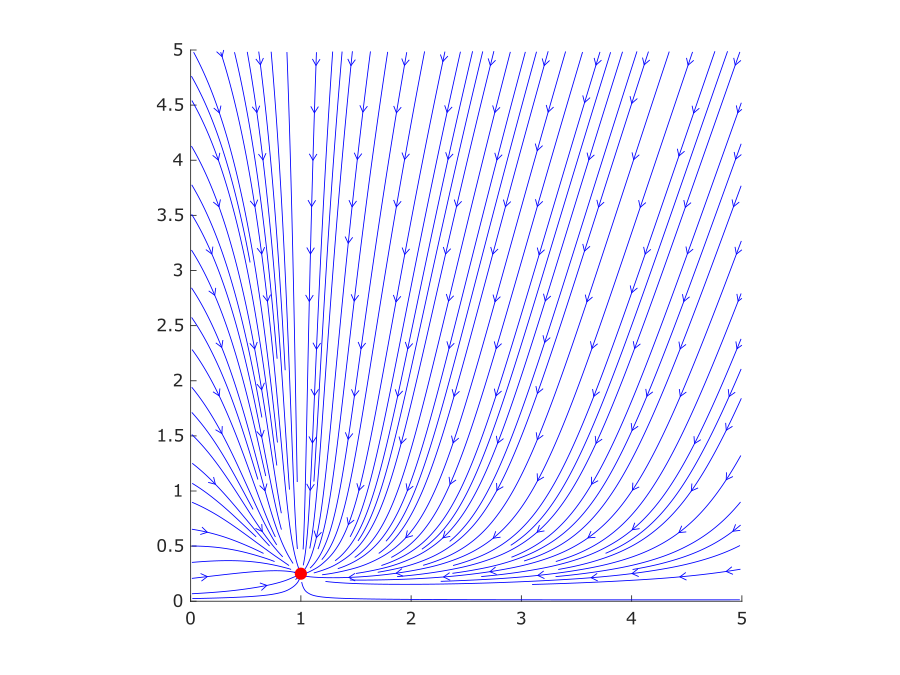}}\hfill
    \subfloat[Full unit]{\includegraphics[width = 0.5\textwidth]{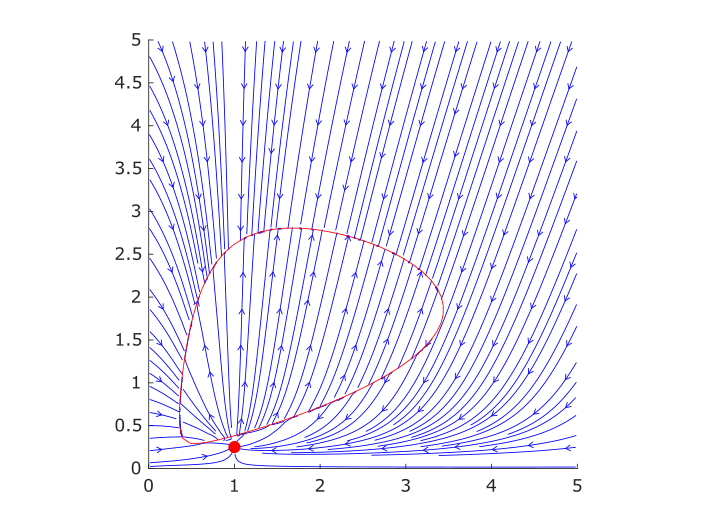}}
    \caption{Phase portraits for the base and full unit.}
    \label{Phase Portraits}
\end{figure}

The curve of infinitely many positive steady states is the solution to the common term(scalar polynomial). The fixed point $(1,0.25)$ is the solution to the set vector polynomials. The phase portraits of the full system given in Figure~\ref{Phase Portraits} show a half-stable limit cycle, with a stable fixed point $(1, 0.25)$ just outside the limit cycle. The phase portrait of the base unit also has a single stable fixed point at $(1, 0.25)$

%\subsubsection{Network Analysis}

We now describe the properties of the network corresponding to the full unit. The network is not weakly reversible, has a deficiency of eleven and is endotactic, but not strongly endotactic. Further, the network possesses critical siphons implying that the dynamics generated by it is not structurally persistent. 

In addition, there exists no weakly reversible mass-action system that is dynamically equivalent to the mass-action system generated by the full unit. For a proof, refer to ~\ref{proof}.\ref{example_1} in the Appendix.

\subsection{Example 2}\label{Exp2}
\begin{figure}[h!]
    \centering
    \includegraphics[scale=0.3]{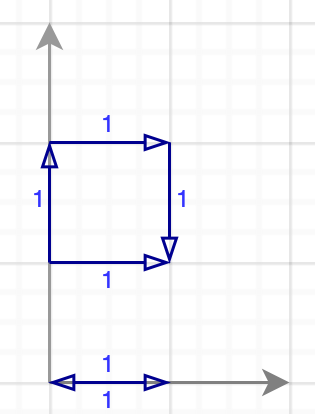} 
    \caption{Base Unit 2}
    \label{Base Unit2}
\end{figure} 
\begin{figure}
    \centering
    \includegraphics[width = 0.23\textwidth, height = 150pt]{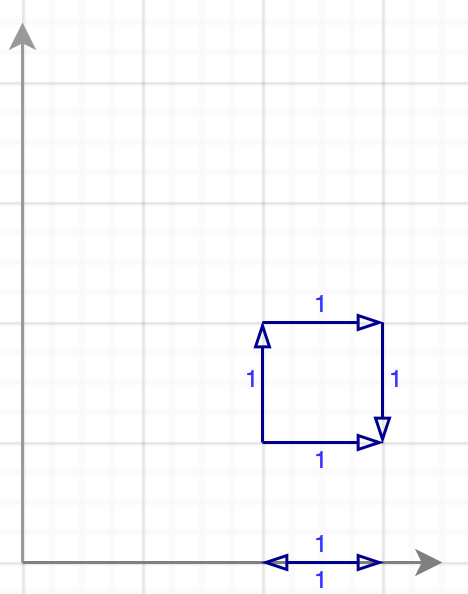}\hfill
    \includegraphics[width = 0.23\textwidth, height = 150pt]{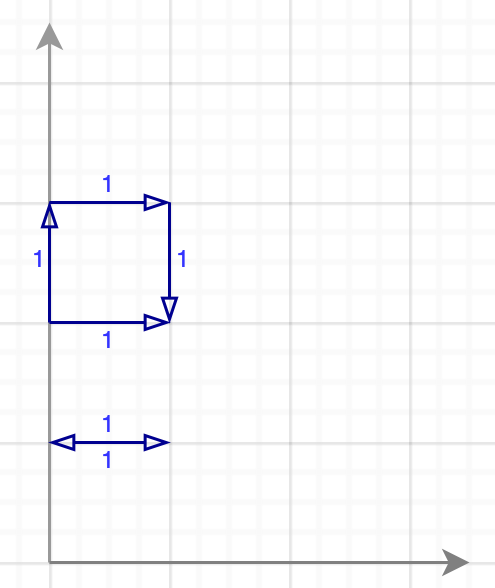}\hfill
    \includegraphics[width = 0.23\textwidth, height = 150pt]{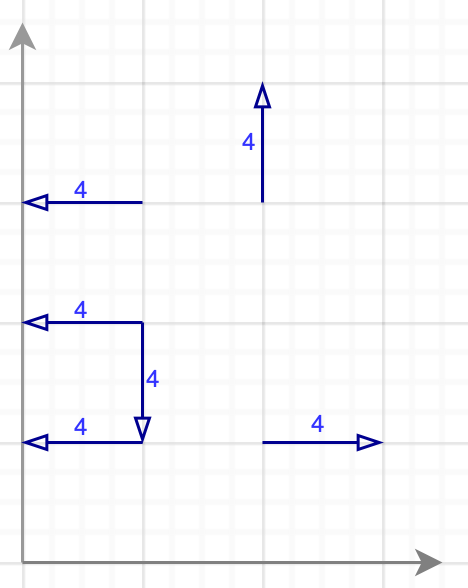}\hfill
    \includegraphics[width = 0.23\textwidth, height = 150pt]{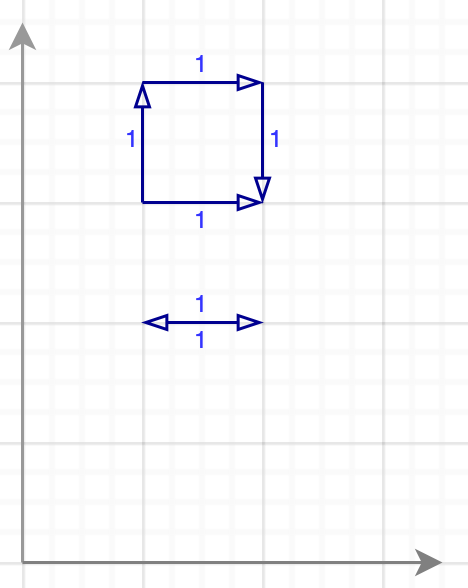}\hfill
    \caption{The different transformations of the base unit given in Figure~\ref{Base Unit2} that combine to form the endotactic network with infinite number of positive steady states}
    \label{Transformations2}
\end{figure}

We start with a base unit given in Figure~\ref{Base Unit2}, where the rate constants corresponding to all reactions are all set to 1. This base unit is endotactic but not weakly reversible. The base unit is not weakly reversible, has a deficiency of two and is endotactic, but not strongly endotactic. Further, the network possesses critical siphons implying that the dynamics generated by it is not structurally persistent. 
The properties mentioned above can be verified with the CoNtRol software~\cite{donnell2014control}.  

\begin{figure}[h!]
    \centering
    \includegraphics[width=0.3\textwidth]{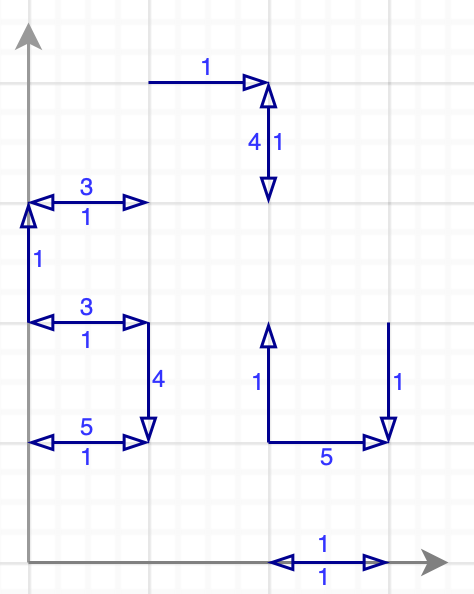}
    \caption{Full unit for Example 2 is given by combining all the transformations in Figure~\ref{Transformations2} and later simplifying it.}
    \label{fullSystem2}
\end{figure}

Further there exists no weakly reversible mass-action system that is dynamically equivalent to the mass-action system generated by the base unit. This base unit is then modified by using the transformations given in Figure~\ref{Transformations2}. Combining these transformations we obtain an endotactic mass action system represented in Figure~\ref{fullSystem2}. Its dynamics is given by:

\begin{equation*}
    \begin{aligned}
        \dot{x}&=(x^2+xy^2+y-4xy)(1-x+y+y^2) \\
        \dot{y}&=(x^2+xy^2+y-4xy)(y-xy^2)
    \end{aligned}
\end{equation*}

\subsubsection{Phase Plane Analysis}
The positive steady states of the dynamical system are given by:
\begin{equation*}
    \{(x,y)^T \in \mathbb{R}_+^2:x^2+xy+y-4xy = 0 \}  \cup \{(1.839, 0.544)^T\}
\end{equation*}
\begin{figure}
    \centering
    \subfloat[Base Unit]{\includegraphics[width = 0.5\textwidth]{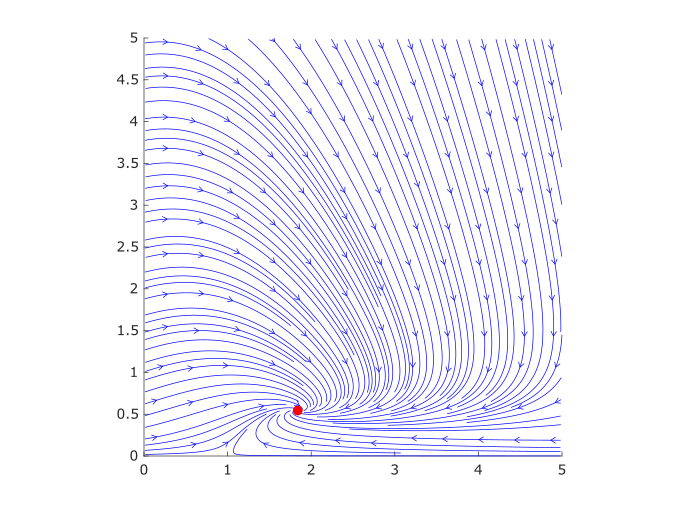}}\hfill
    \subfloat[Full unit]{\includegraphics[width = 0.5\textwidth]{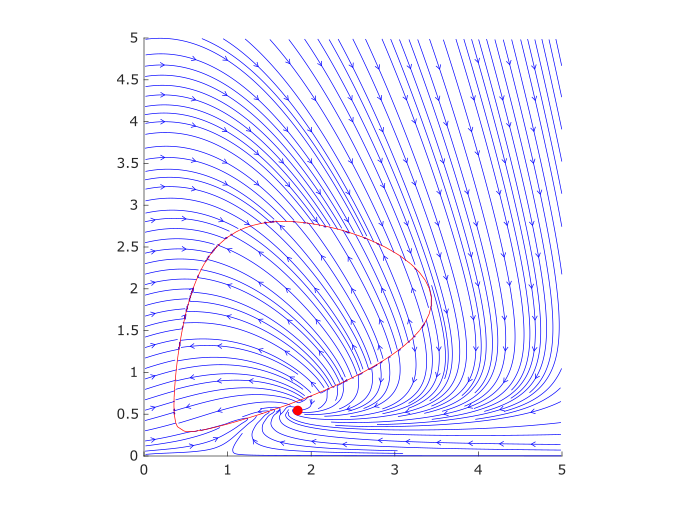}}
    \caption{Phase portraits for base and full unit.}
    \label{Phase Portraits2}
\end{figure}

The curve of the infinitely many steady states is the solution to the common term(scalar polynomial) and the point $(1.839, 0.544)$ is the solution of the set of vector polynomials. The phase portrait of the full system given in Figure~\ref{Phase Portraits2} shows a stable fixed point $(1.839, 0.544)$ outside the curve of fixed points(which is a semi-stable limit cycle). The phase portrait of the base unit has a single stable fixed point at $(1.839, 0.544)$.

%\subsubsection{Network Analysis}

We now describe the properties of the network corresponding to the full unit. The network is not weakly reversible, has a deficiency of nine and is endotactic, but not strongly endotactic. Further, the network possesses critical siphons implying that the dynamics generated by it is not structurally persistent. 

Further, there exists no weakly reversible mass-action system that is dynamically equivalent to the mass-action system generated by the full unit. For a proof, refer~\ref{proof}.\ref{example_2} in the Appendix.

\subsection{Example 3}

\begin{figure}[h!]
    \centering
    \includegraphics[scale=0.4]{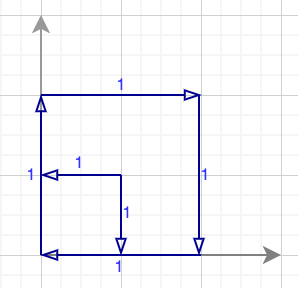}
    \caption{Base Unit}
    \label{Base Unit 3}
\end{figure} 

\begin{figure}
    \centering
    \includegraphics[width = 0.23\textwidth, height = 150pt]{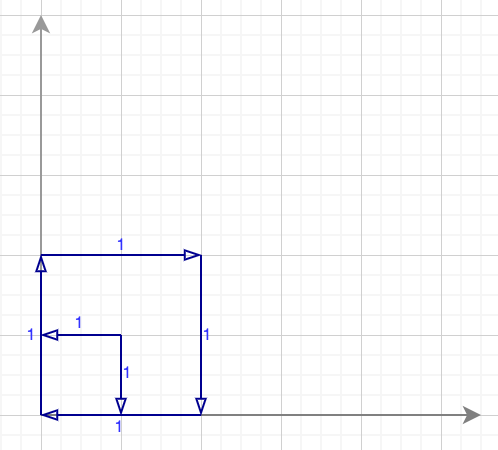}\hfill
    \includegraphics[width = 0.23\textwidth, height = 150pt]{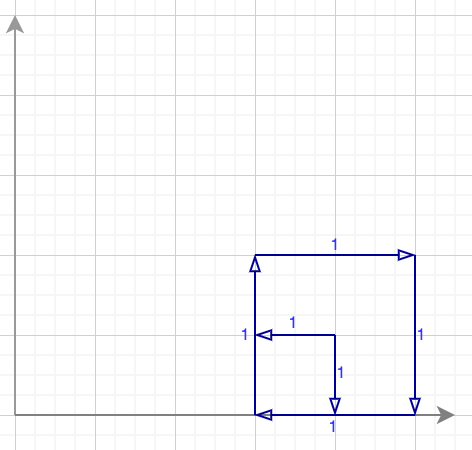}\hfill
    \includegraphics[width = 0.23\textwidth, height = 150pt]{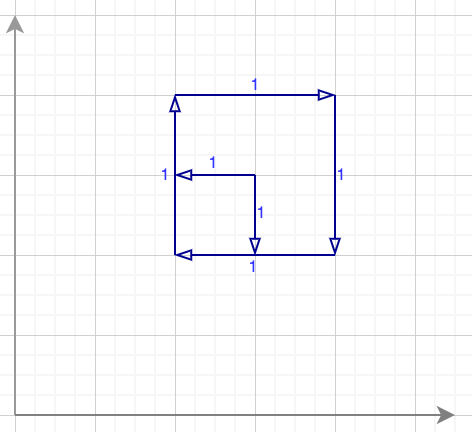}\hfill
    \includegraphics[width = 0.23\textwidth, height = 150pt]{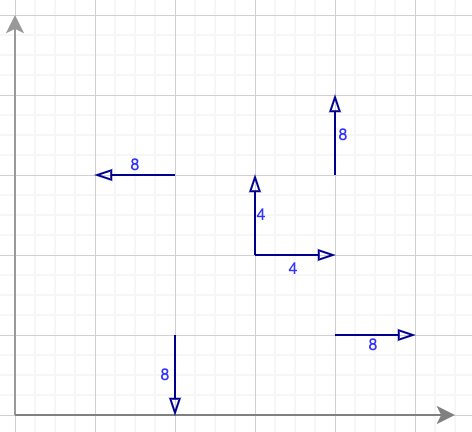}\hfill
    \caption{The different transformations of the base unit given in Figure~\ref{Base Unit 3} that combine to form the endotactic network with infinite number of positive steady states}
    \label{Transformations3}
\end{figure}

\begin{figure}[h!]
    \centering
    \includegraphics[width=0.3\textwidth]{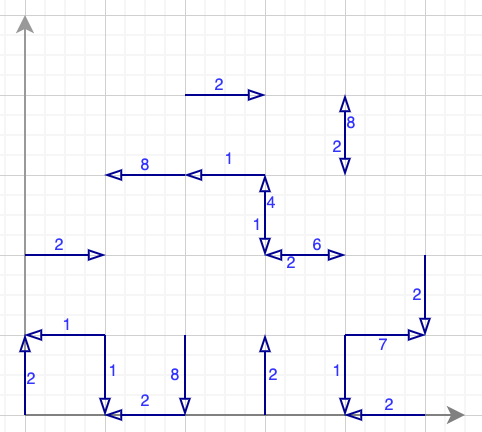}
    \caption{Full unit for Example 3 is given by combining all the transformations given in Figure~\ref{Transformations3} of the base unit and simplifying it.}
    \label{fullSystem3}
\end{figure}

In this example, we start with a base unit given in Figure~\ref{Base Unit 3}, where the rate constants corresponding to all reactions are set to 1. This base unit is strongly endotactic but not weakly reversible. The network has a deficiency of four and is strongly endotactic. Further, the network possesses no critical siphons implying that the dynamics generated by it is structurally persistent. The properties mentioned above can be verified with the CoNtRol software~\cite{donnell2014control}.

Further there exists no weakly reversible mass-action system that is dynamically equivalent to the mass-action system generated by the base unit. This base unit is then modified by using the transformations given in Figure~\ref{Transformations3}. Combining these transformations we obtain a strongly endotactic mass action system depicted in Figure~\ref{fullSystem3} The dynamics of the full unit is given by:

\begin{equation*}
    \begin{aligned}
        \dot{x}&=(1+x^3+x^2y^2-4x^2y)(2y^2-2x^2-xy)\\
        \dot{y}&=(1+x^3+x^2y^2-4x^2y)(2-2x^2y^2-xy)
    \end{aligned}
\end{equation*}

\subsubsection{Phase Plane Analysis}
The steady states of the system is given by: 
\begin{equation*}
    \{(x,y)^T \in \mathbb{R}_+^2: 1+x^3+x^2y^2-4xy = 0 \}  \cup \{(0.781, 1)^T\}
\end{equation*}

\begin{figure}
    \centering
    \subfloat[Base Unit]{\includegraphics[width = 0.5\textwidth]{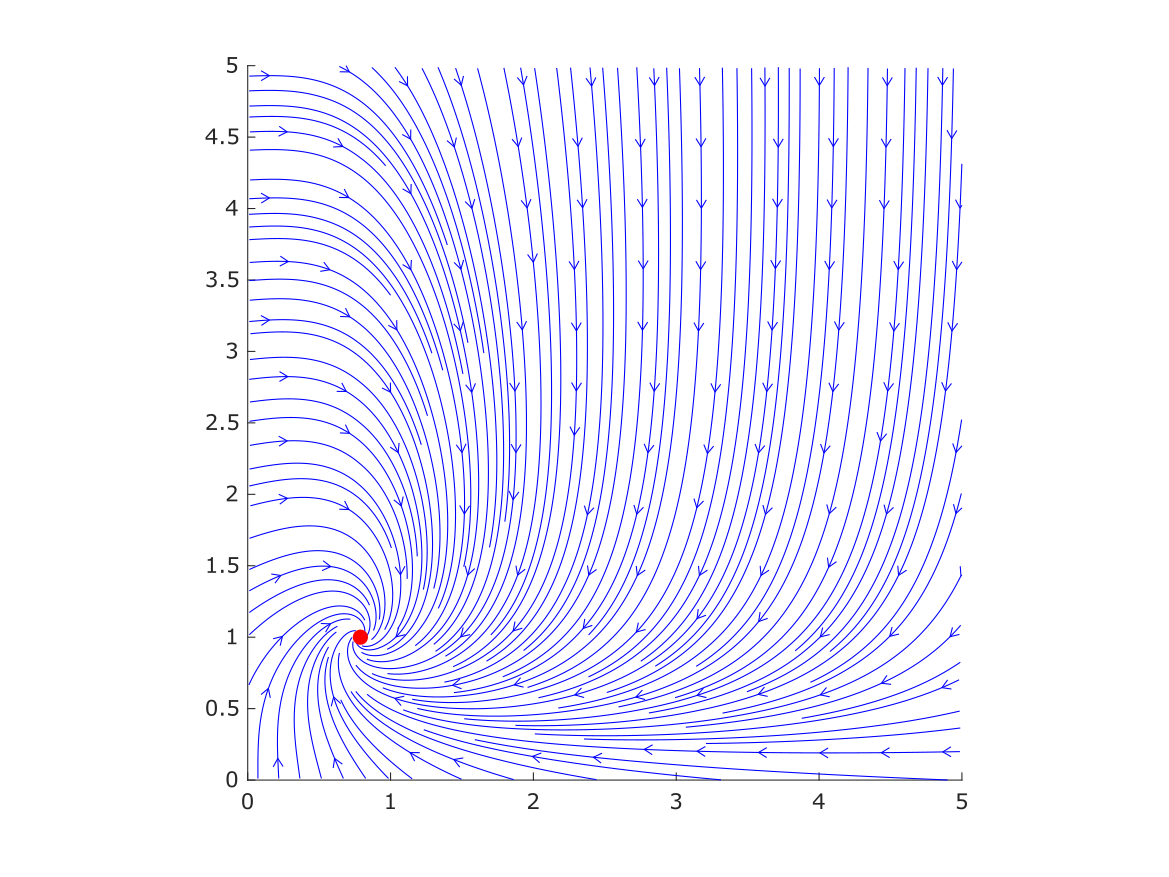}}\hfill
    \subfloat[Full unit]{\includegraphics[width = 0.5\textwidth]{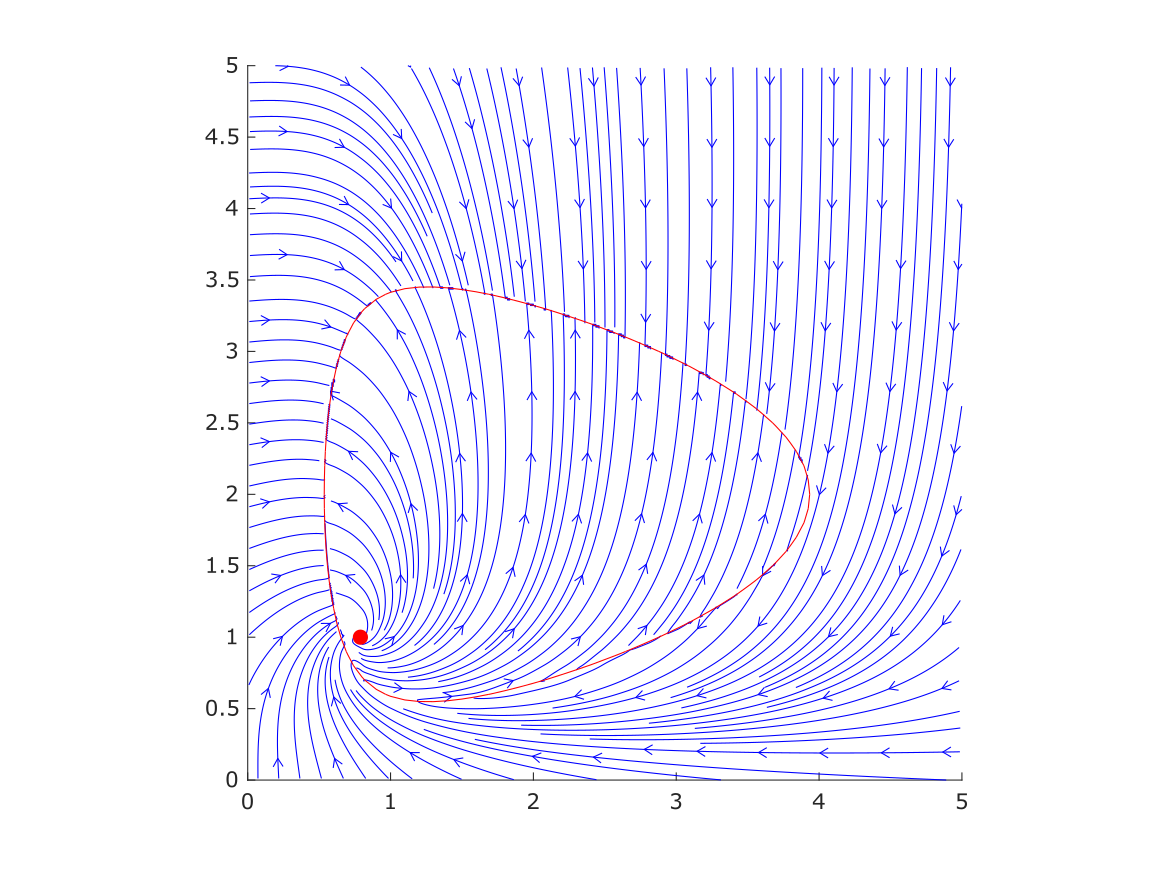}}
    \caption{Phase portraits for the base and full unit}
    \label{Phase Portraits3}
\end{figure}

The curve of infinite steady states is the solution to the common term(scalar polynomial) and the point $(0.781, 1)$ is the solution for the set of vector polynomials. The phase portrait of the full system is given in Figure \ref{Phase Portraits3} show a unstable fixed point, surrounded by a stable limit cycle(which is a curve of infinite steady states.). The phase portrait of the base unit dynamical system have a single stable fixed point at (0.781, 1).

%\subsubsection{Network Analysis}

We now describe the properties of the network corresponding to the full unit. The network is not weakly reversible, has a deficiency of fifteen and is strongly endotactic. Further, the network possesses no critical siphons implying that the dynamics generated by it is structurally persistent.

\section{Networks in three dimensions}\label{sec:higherDimensions}

The examples that we have analyzed so far are two dimensional. We now construct a three dimensional endotactic reaction network that is not weakly reversible but possesses a curve of infinitely many positive steady states.

\subsection{Example}

\begin{figure}[h!]
    \centering
    \includegraphics[scale=0.3]{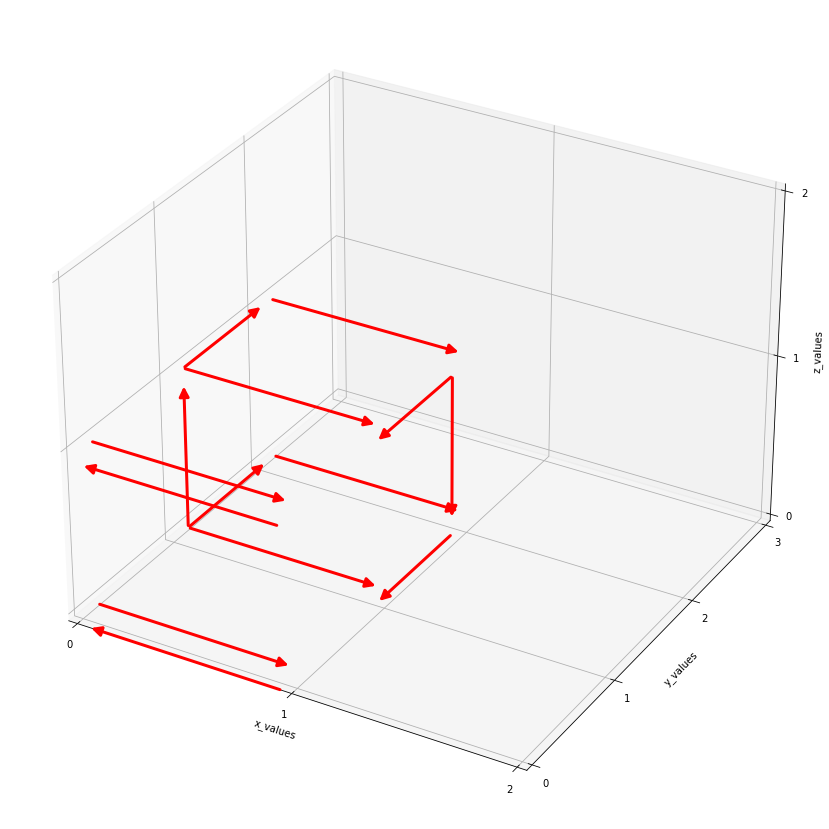}
    \caption{Three Dimensional Base Unit}
    \label{3dBaseUnit}
\end{figure}

To come up with a base unit in three dimensions, we start with the two dimensional base unit given in Figure~\ref{Base Unit2} and then combine it with a copy of itself transformed along the z-axis. To this system we add two reactions to connect the two planar vertices so that the reaction network looks like Figure~\ref{3dBaseUnit}. 

The resultant network is not weakly reversible, has a deficiency of six and is endotactic, but not strongly endotactic. Further, the network possesses critical siphons implying that the dynamics generated by it is not structurally persistent. The properties mentioned above can be verified with the CoNtRol software.\cite{donnell2014control}. \\

Note that there exists no weakly reversible mass-action system that is dynamically equivalent to the mass-action system generated by the base unit. The dynamics corresponding to the base unit is given by:
\begin{equation*}
    \begin{aligned}
        \dot{x}&=(1-x+y+y^2+z-zx+zy+zy^2) \\
        \dot{y}&=(y-xy^2+zy-zxy^2) \\
        \dot{z}&=(y-xy^2z)
    \end{aligned}
\end{equation*}

\begin{figure}
    \centering
    \includegraphics[width = 0.3\textwidth, height = 150pt]{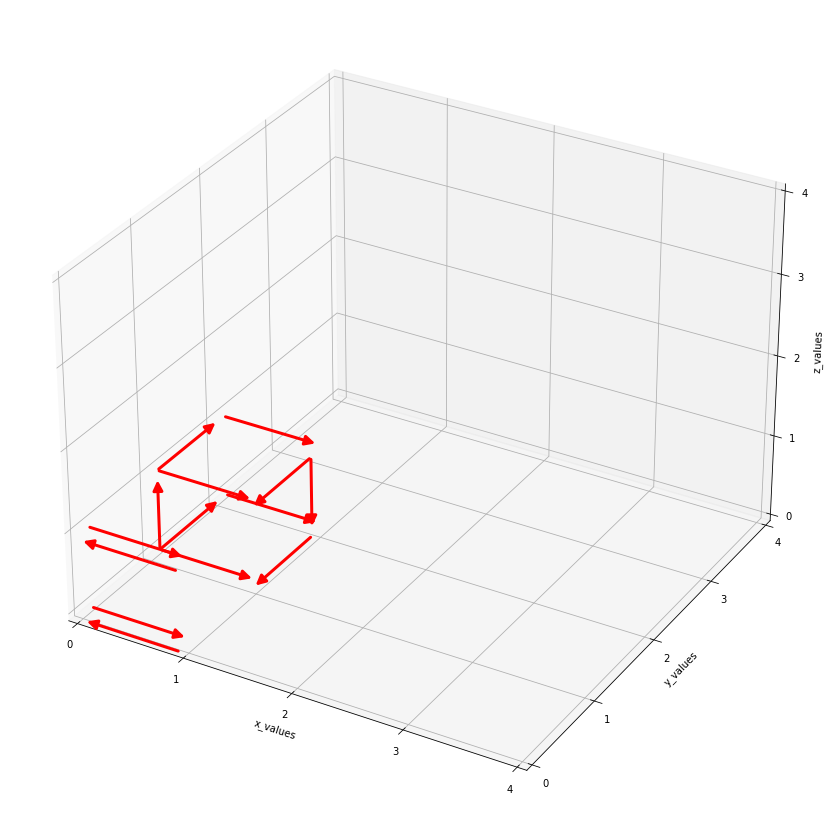}\hfill
    \includegraphics[width = 0.3\textwidth, height = 150pt]{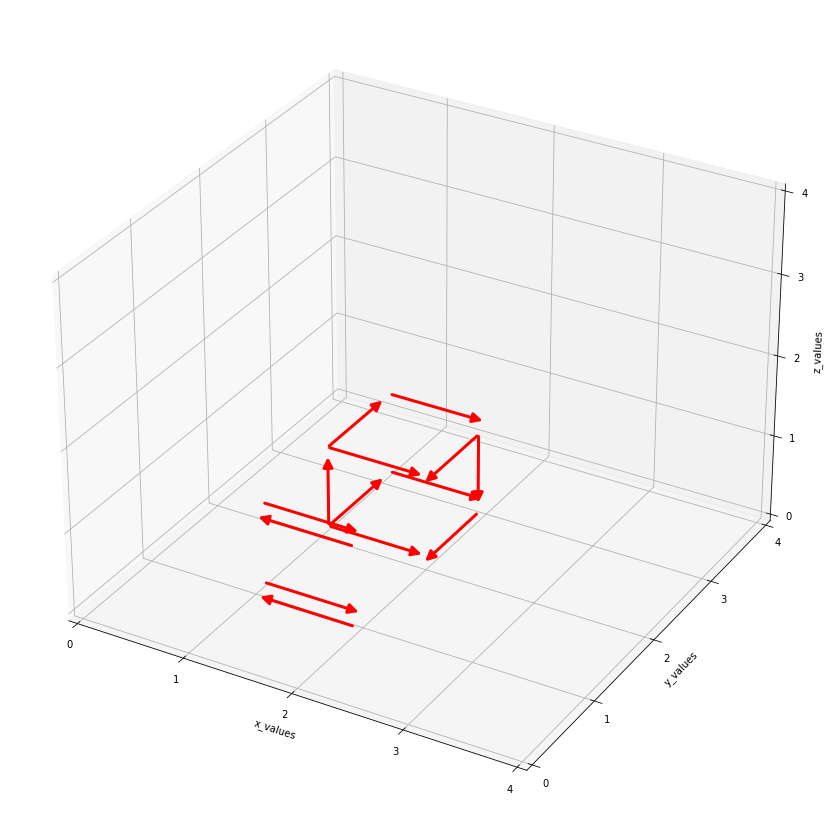}\hfill
    \includegraphics[width = 0.3\textwidth, height = 150pt]{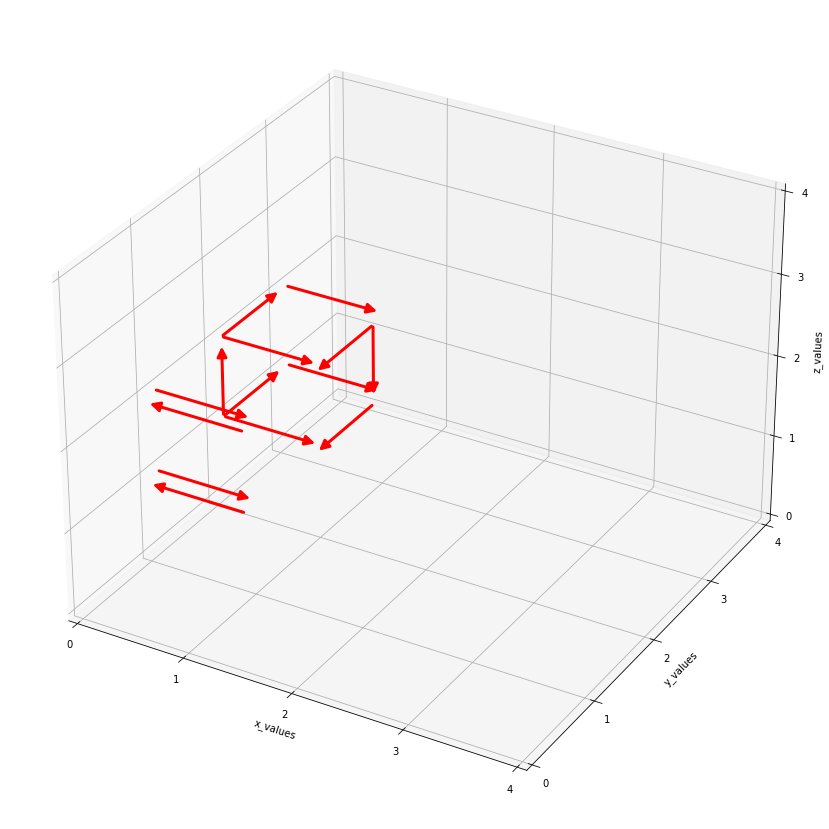}\hfill \\
    \includegraphics[width = 0.3\textwidth, height = 150pt]{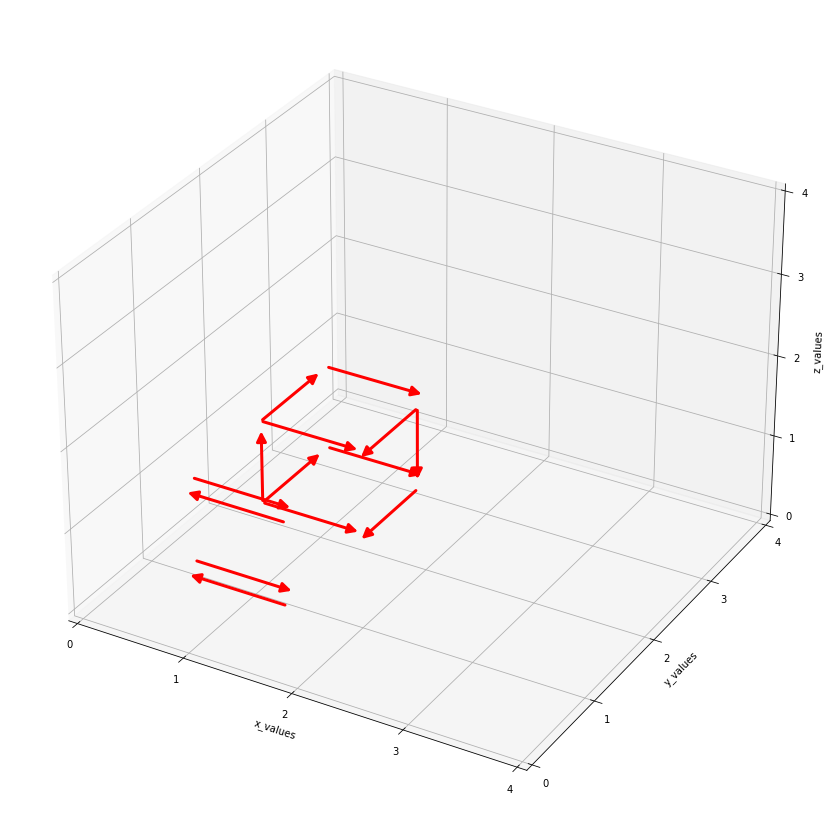}\hfill
    \includegraphics[width = 0.3\textwidth, height = 150pt]{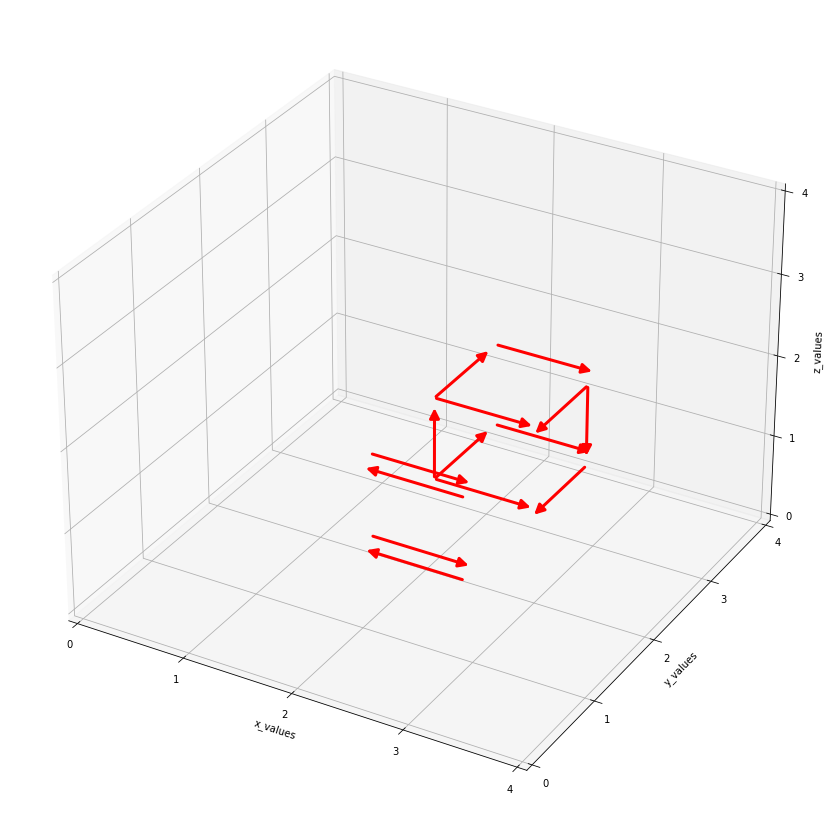}\hfill
   \includegraphics[width = 0.3\textwidth, height = 150pt]{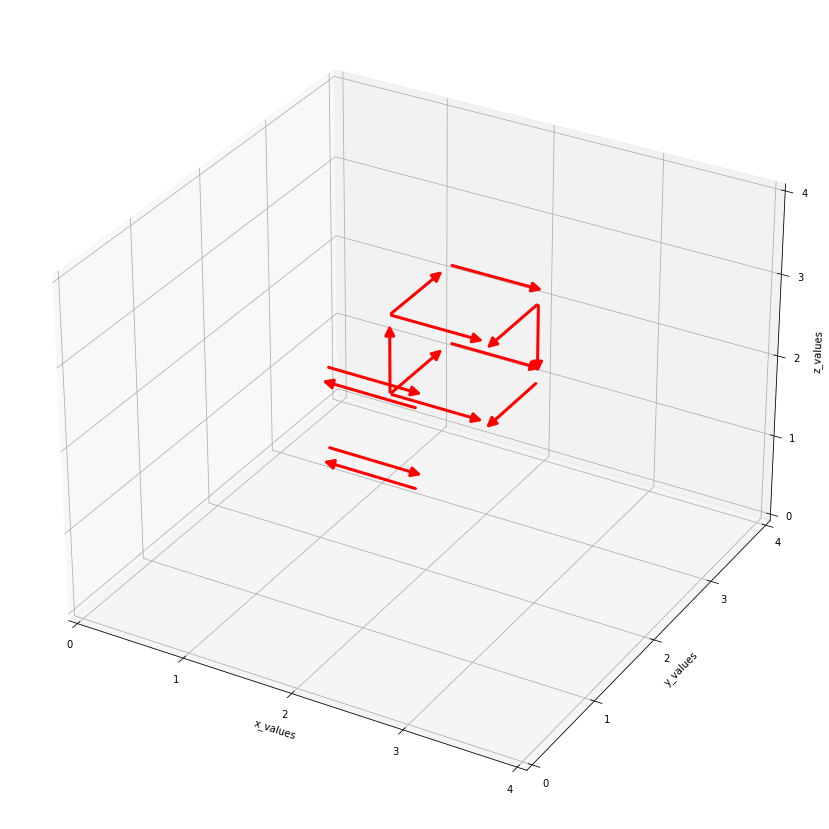}\hfill \\
    \includegraphics[width = 0.3\textwidth, height = 150pt]{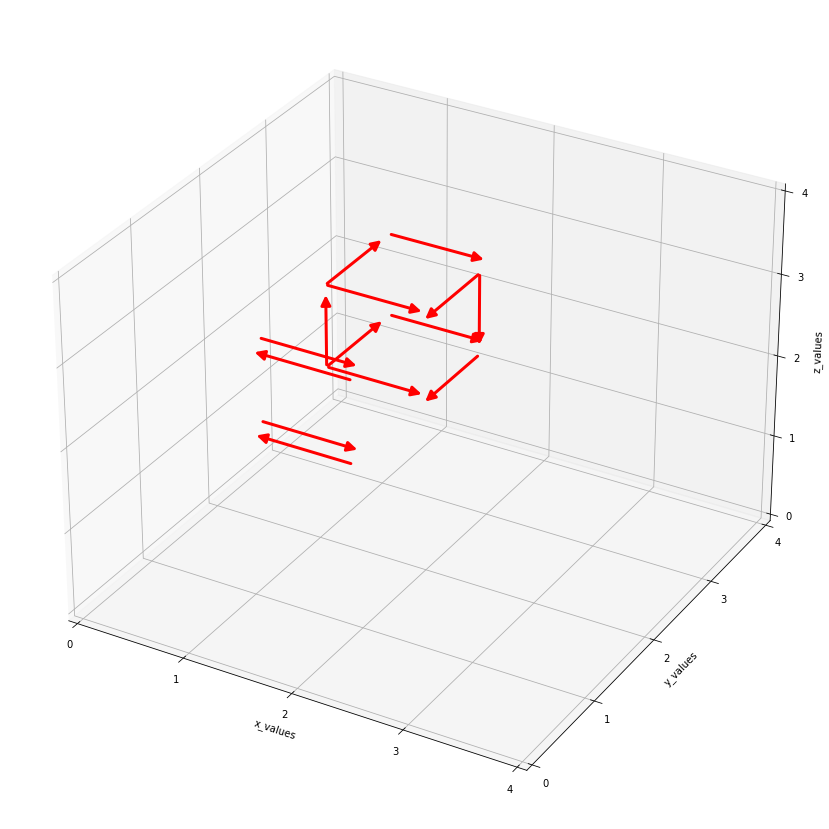}\hfill
    \includegraphics[width = 0.3\textwidth, height = 150pt]{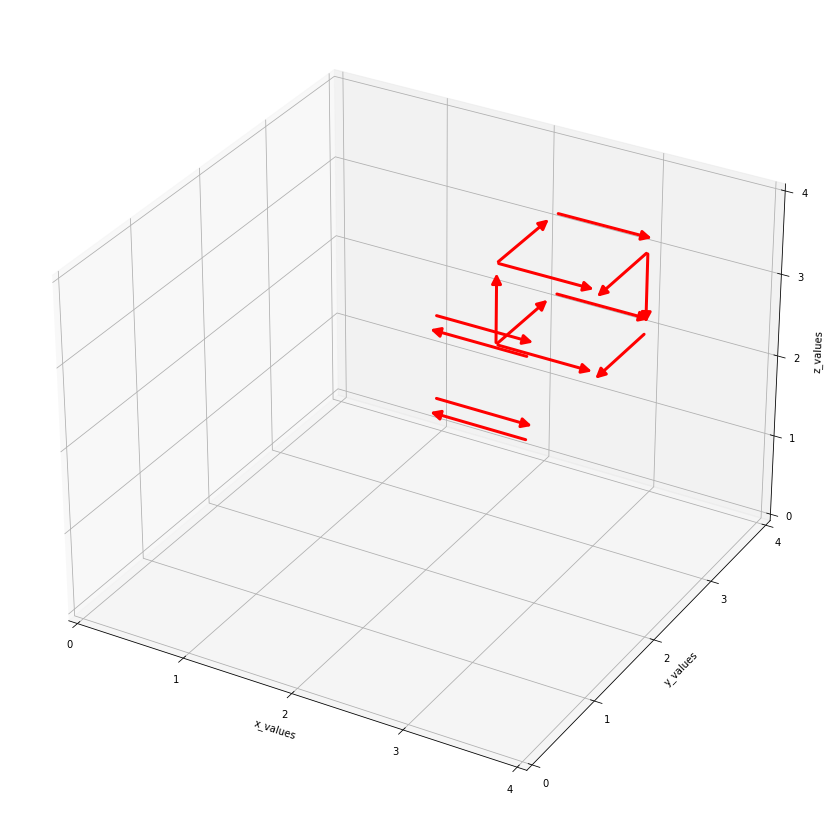}\hfill
    \includegraphics[width = 0.3\textwidth, height = 150pt]{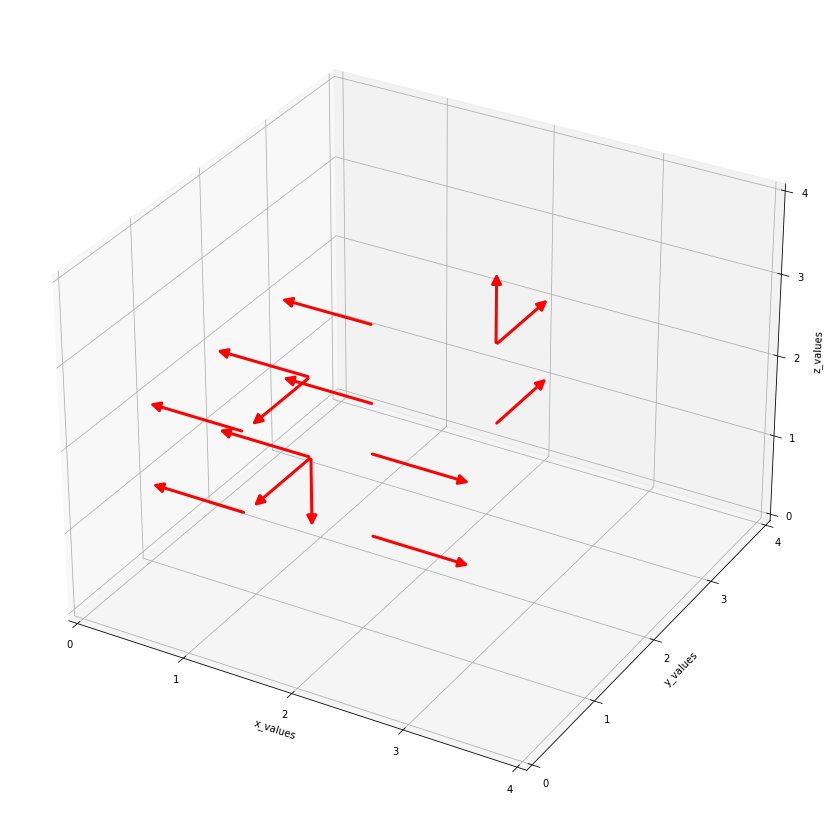}\hfill
    \caption{The different transformations of the base unit given in Figure~\ref{3dBaseUnit} that combine to form the endotactic network with infinite number of positive steady states}
    \label{Transformations3d}
\end{figure}

The base unit is then modified using certain transformations given in \ref{Transformations3d}. The scalar polynomial for these transformations is given by:
\begin{equation*}
    h(x,y,z) = (1+xy+yz+xz+x^2yz+xy^2z+xyz^2+x^2y^2z^2-15xyz)
\end{equation*}
Combining these transformations, we obtain the mass-action system depicted in figure \ref{fig:3dfinalUnit}.

\subsubsection{Phase Plane Analysis}

The steady states of the base system are given by:
\begin{equation*}
    \{(1,0,z)^T  \cup {(1.83,0.54,1)^T}\}
\end{equation*}
where $z \geq 0$. \\

The dynamics generated by the base unit has infinitely many positive steady states which is given by a line of fixed points $(x=1, y=0, z \geq 0 )$. The steady states of the full system is given by the set:
\begin{equation*}
    A=\{(x,y,z)^T \in \mathbb{R}^3_{>0}:\{h(x,y,z) = 0\} \cup (1,0,z)^T\cup {(1.83,0.54,1)^T}\} 
\end{equation*}

\begin{figure}[h!]
    \centering
    \begin{subfigure}[t]{0.5\textwidth}
        \centering
        \includegraphics[trim={0cm 0cm 2.5cm 0cm},clip,width = 0.65\textwidth]{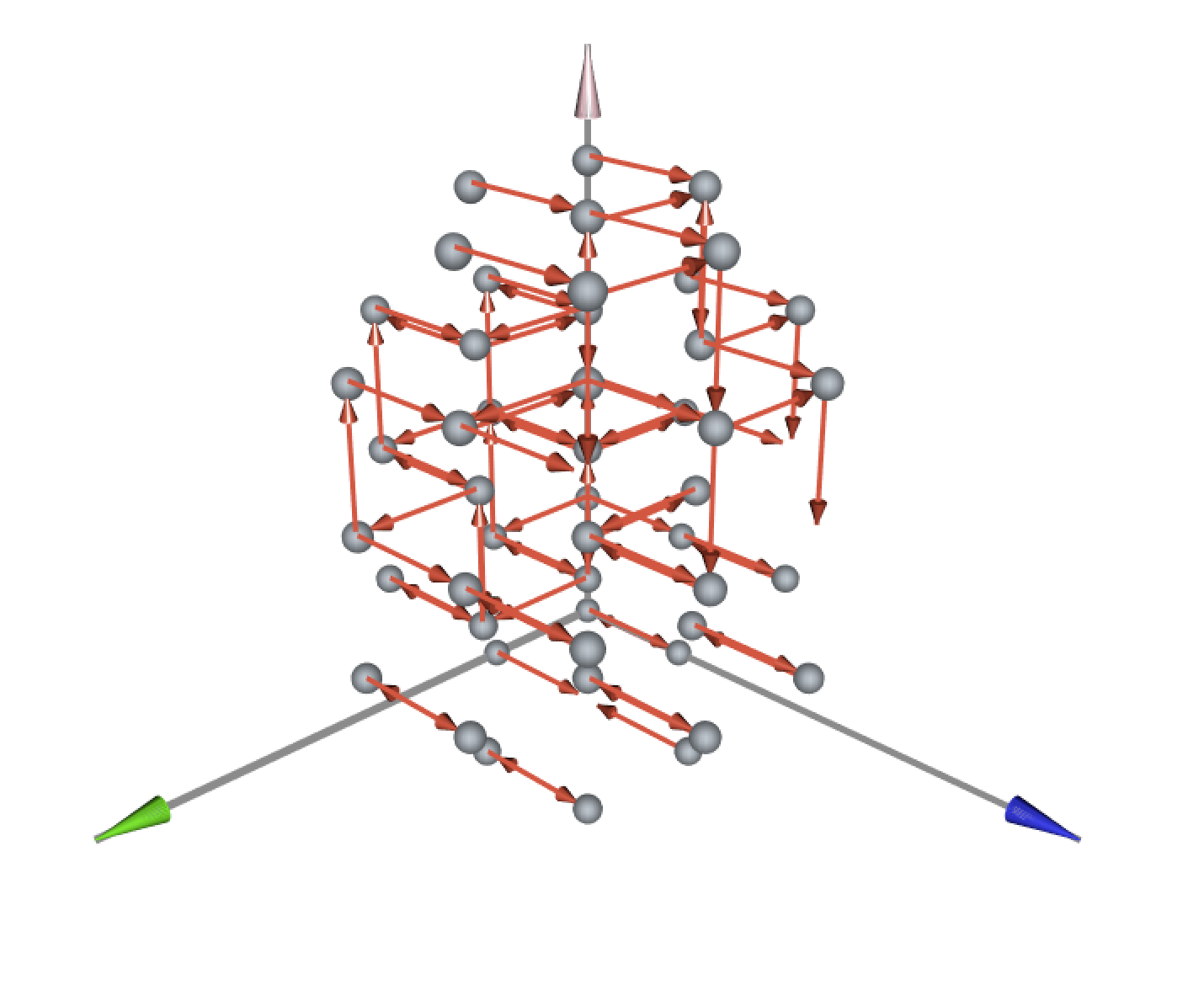}
        \caption{Three dimensional final unit.}
        \label{fig:3dfinalUnit}
    \end{subfigure}%
    \begin{subfigure}[t]{0.5\textwidth}
        \centering
        \includegraphics[trim={5cm 0cm 5cm 0cm},clip,width = 0.65\textwidth]{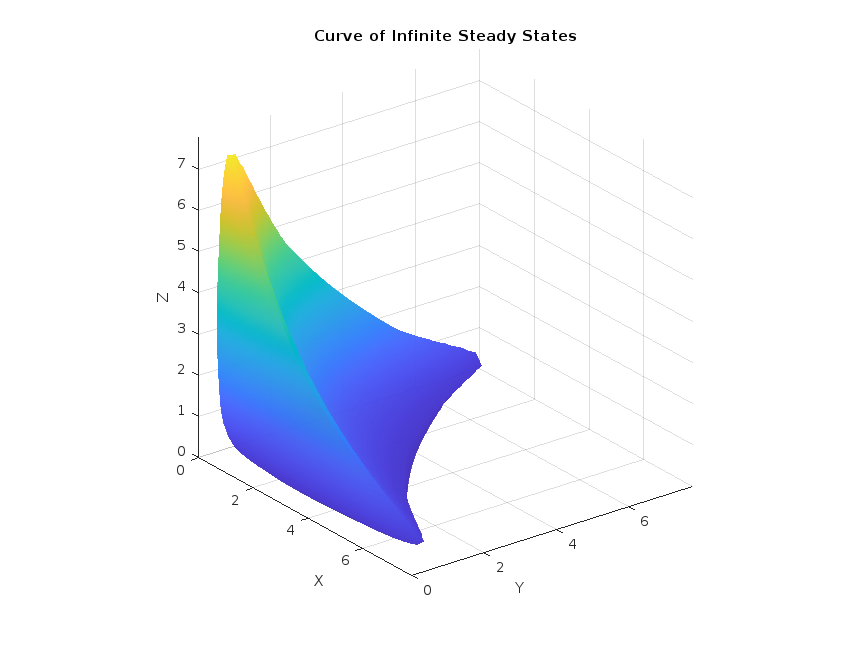}
        \caption{Curve of infinitely many positive steady states.}
        \label{fig:3dCurveofInfiniteSteayStates}
    \end{subfigure}%
    \caption{Three dimensional final unit and the curve of infinitely many positive steady states.}

\end{figure}

The full system gives us a curve of infinitely many steady states in $R^3_{>0}$, which is represented in Figure~\ref{fig:3dCurveofInfiniteSteayStates}. \\

We now describe the properties of the network corresponding to the full unit. The network is not weakly reversible, has a deficiency of forty four and is endotactic, but not strongly endotactic. Further the network possesses critical siphons implying that the dynamics generated by it is not structurally persistent.

%\subsubsection{Network Analysis}

%and an interactive 3D figure is given \href{https://samay-kothari.github.io/ResearchWork/}{here}. 
%The steady states of the full unit system are given by:
%\begin{equation*}
%    \{(x,y,z)^T \in \mathbb{R}_+^3:1+xy+yz+xz+x^2yz+xy^2z+xyz^2+x^2y^2z^2-15xyz = 0 \}  \cup 
%\end{equation*}

\begin{comment}
The final system after combining all the transformations is given \href{https://samay-kothari.github.io/ResearchWork/}{here}, and some snapshots of the same is given in figure \ref{fig:3dsnapshots}. \\ \\
\begin{figure}
    \centering
    \includegraphics[trim={5cm 5cm 5cm 5},clip,width = 0.45\textwidth]{images/example3d/pic00.png}
    \includegraphics[trim={5cm 5cm 5cm 5},clip,width = 0.45\textwidth]{images/example3d/pic02.png} \\
    \includegraphics[trim={5cm 5cm 5cm 5},clip,width = 0.45\textwidth]{images/example3d/pic03.png}
    \includegraphics[trim={5cm 5cm 5cm 5},clip,width = 0.45\textwidth]{images/example3d/pic01.png}
    \caption{Full unit for the 3d example after combining all the transformations given in \ref{Transformations3} of the base unit \ref{3dBaseUnit}}
    \label{fig:3dsnapshots}
\end{figure}

\end{comment}

\section{Appendix}\label{Appendix}
\subsection{Possible transformations of a reaction vector}

In what follows, we list some transformations that can be done on a reaction vector and its rate constant so that the resultant system is dynamically equivalent to the original system.

\subsubsection{Length transformation}\label{scalarTransformation}

\begin{figure}[h!]

    \begin{subfigure}[t]{0.49\textwidth}
        \centering
        \includegraphics[width=0.5\textwidth]{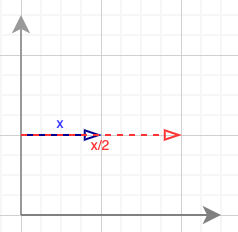}
        \caption{Length transformations}
        \label{fig:lengthTransformation}
    \end{subfigure}
    \begin{subfigure}[t]{0.49\textwidth}
        \centering
        \includegraphics[width=0.5\textwidth]{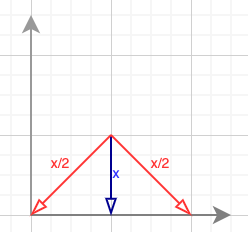}
        \caption{Diagonal decomposition of reaction vector}
        \label{fig:diagonalDecomposition}
    \end{subfigure}%
    \caption{Transformations of the reaction vectors}
    \label{fig:vectorTransformations}
\end{figure}

This kind of transformation refers to the case where we modify the length of the reaction vector and its rate constant suitably, whilst keeping the direction of the reaction vector the same to ensure dynamical equivalence. Figure~\ref{fig:lengthTransformation} gives an example of a scalar transformation of a reaction vector.

\subsubsection{Diagonal decomposition}\label{diagonalDecomposition}

In this case, we take the reaction vector and decompose it into two reaction vectors, such that the original reaction vector bisects the two new vector. The rate constants of the new vectors are modified accordingly. Figure~\ref{fig:diagonalDecomposition} gives an example of the this vector transformation.

\subsection{Proof that there exists no weakly reversible mass-action systems that are dynamically equivalent to the mass-action systems generated by the full units in Examples 1 and 2}\label{proof}
% \begin{figure}
%     \centering
%     \includegraphics[width = 0.5 \textwidth]{images/appendix/FinalSystem1.drawio.png}
%     \caption{Example 1 non dynamical equivalence to weakly reversible}
%     \label{appendixExample1}
% \end{figure}
%ranslatio
\begin{figure}
    \centering
    \begin{subfigure}[t]{0.45\textwidth}
        \centering
        \includegraphics[width=0.6\textwidth]{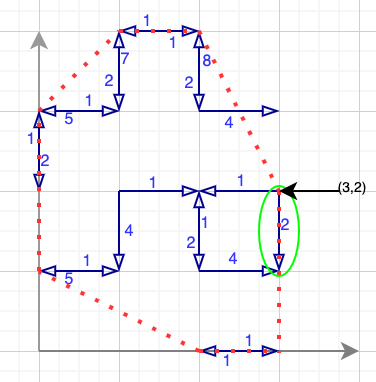}
        \caption{Full unit of example 1}
        \label{fig:appendixExample1}
    \end{subfigure}%
    \begin{subfigure}[t]{0.45\textwidth}
        \centering
        \includegraphics[width=0.6\textwidth]{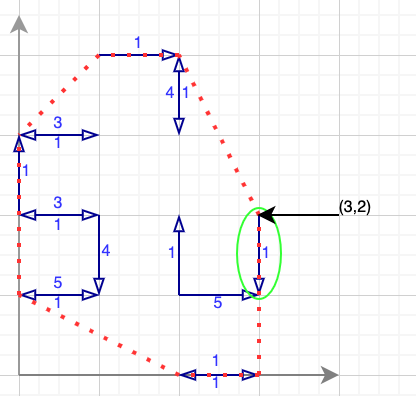}
        \caption{Full unit of example 2}
        \label{fig:appedixExample2}
    \end{subfigure}%
    \caption{Full systems for example 1 and 2. Both these systems are endotactic, but are neither weakly reversible nor strongly endotactic. The boundary of the convex hull is marked with a dotted line. Further, there exists no weakly reversible mass-action systems that are dynamically equivalent to the mass-action systems generated by the these networks.}
    \label{fig:appendixExamples}
\end{figure}

\begin{enumerate}

\item\label{example_1} \textbf{Example 1:} In Figure~\ref{fig:appendixExample1} we have the full unit given in Example 1. The boundary of the convex hull of the full unit is marked with a dotted red line. The problematic reaction is circled in green. There is no way which to return to the starting vertex, i.e. the point $(3,2)$. 

If we try to decompose the green reaction according to~\ref{diagonalDecomposition}, the new vector that point outwards of the convex hull(on the right side), would again lead us to a vertex from where we cannot return.

If we try extending the reaction as we did in~\ref{scalarTransformation}, the new reaction vector takes us to the point $(3,0)$. From  the point $(3,0)$, we cannot go back to the upper part of the reaction network due to the reversible reaction vector between $(2,0)$ and $(3,0)$ 
Therefore, there exists no weakly reversible dynamical system that is dynamically equivalent to the mass-action system generated by the final system in Example 1.

\item \label{example_2} \textbf{Example 2:} In Figure~\ref{fig:appedixExample2} we have the full unit given in Example 2. The boundary of the convex hull of the full unit is marked with a dotted red line. The problematic reaction is circled in green. Arguing as in Example 1, there is no way to return to the starting vertex, i.e. the point $(3,2)$. Therefore, there exists no weakly reversible dynamical system that is dynamically equivalent to the the mass-action system generated by the final system in Example 2.

\end{enumerate}

\section{Discussion}\label{sec:discussion}

% \section{Minutes of the Meeting}
% \begin{enumerate}
%     \item Add the three endotactic networks documented, and look out for more, and the way they are constructed.
%     \item Can a weakly reversible mass-action system have infinitely many positive steady states without having a common factor on the right-hand side of its differential equations
% \end{enumerate}
% \section{Next Project}
% \begin{itemize}
%     \item Finding relationship between conservative and consistent reaction networks.
%     \item Study about petrinets (special notion of conservative petrinets)
% \end{itemize}
% \subsection{Important Points}
% \begin{itemize}
%     \item Endotactic networks are consistent
%     \item Consistent networks with zero deficiency are weakly reversible and thus endotactic as well.
%     \item \textbf{Definition:} Reaction network is consistent, if there is a positive combination(all coefficients are $>$ 0) of reaction networks that sum to 0.
%     Simplest: x+y $\rightarrow$ 0 and x+y $\rightarrow$ 2x+2y \\
%     \item \textbf{Definition:} Reaction networks in which we have at least one positive vector that is orthogonal to the stoichiometric subspace. Simplest: x+y $\rightarrow$ 2y
% \end{itemize}

% \begin{align}z
% \begin{split}
% X + Y + Z &\rightarrow 3Z \\
% 2X + Y &\rightarrow 4Z
% \end{split}
% \end{align}

% The dynamical system it generates is given by

% \begin{align}
% \begin{split}
% \dot{x} &= -xyz -2x^2y \\
% \dot{y} &= -xyz - x^2y \\
% \dot{z} &= 2xyz  + 4x^2y
% \end{split}
% \end{align}

The study of steady states is crucial to the analysis of dynamical systems. In particular, questions about the existence and the number of steady states are of special importance. In particular, Boros~\cite{boros2019existence} has established that for weakly reversible networks, there exists a positive steady state in every stoichiometric compatibility class. Gopalkrishnan, Miller ans Shiu have shown that strongly endotactic networks also possess at least positive steady state in every stoichiometric compatibility class. Boros, Craciun, and Yu~\cite{boros2020weakly}  have shown that weakly reversible mass-action systems can possess infinitely many positive steady states. We extend this further to show that there exists endotactic and strongly endotactic networks that are not weakly reversible and possess infinitely many positive steady states. In addition, for the base and full units in Examples 1 and 2, there exists no weakly reversible systems that are dynamically equivalent to the mass-action systems generated by them.

This questions raises several intriguing questions for future work. For example,
\begin{enumerate}

\item What is the \emph{minimal} endotactic/strongly endotactic network that possesses infinitely many positive steady states? Here the phrase minimal could mean for e.g. the minimum number of reactions, or the minimum number of complexes or the minimum deficiency.

\item Can one find conditions on the base unit, final unit and the scalar polynomial that ensures that there exists no weakly reversible systems that are dynamically equivalent to the mass-action systems generated by them?

\item Can one find endotactic/strongly endotactic mass-action systems in higher dimensions (dimension greater than four) that are not weakly reversible and possess infinitely many positive steady states ?

\item Is there a  systematic way to construct such networks (i.e., can one design novel algorithms that can help us find such networks in practice?)

\end{enumerate}

\bibliographystyle{amsplain}
\bibliography{Bibliography}

\end{document}